\documentclass[11pt,a4paper]{article}
\usepackage{graphicx}         
\usepackage{cite,url,algorithm}    
\usepackage{amsmath,amsfonts,amssymb}
\usepackage{graphicx,url,algorithm}
\usepackage{color}
\usepackage{algorithm}
\usepackage[noend]{algpseudocode}
\usepackage{pgfplots}
\pgfplotsset{compat=newest}
\pgfplotsset{plot coordinates/math parser=false}
\usepackage{forest}
\usepackage{tikz}
\usetikzlibrary{positioning}
\usepackage{etoolbox}

\newcommand{\proof}{\emph{Proof. }}
\newcommand{\QED}{\hfill\mbox{\rule[0pt]{1.5ex}{1.5ex}}\vspace{0em}\par\noindent}

\usepackage{enumitem}
\renewcommand{\theenumi}{\arabic{enumi}}

\renewcommand{\theenumii}{\arabic{enumii}}

\renewcommand{\theenumiii}{\arabic{enumiii}}

\renewcommand{\theenumiv}{\arabic{enumiv}}

\definecolor{black}{rgb}{0,0,0}
\newcommand{\beqarno}{\begin{eqnarray*}}
\newcommand{\eeqarno}{\end{eqnarray*}}
\newcommand{\st}{\mathrm{s.t.}}
\newcommand{\eqdef}{\triangleq}
\newcommand{\rr}{\mathbb R}
\newcommand{\N}{\mathbb N}

\newcommand{\smallmat}[1]{\left[ \begin{smallmatrix}#1 \end{smallmatrix} \right]}
\newcommand{\card}{\mathop{\rm card}\nolimits}
\newcommand{\Aa}{\mathcal{A}}
\newcommand{\BB}{b}
\newcommand{\DD}{d}
\newcommand{\II}{\mathcal{I}}
\newcommand{\JJ}{\mathcal{J}}
\newcommand{\EE}{\mathcal{E}}
\newcommand{\UU}{\mathcal{U}}

\newcommand{\NN}{\mathcal{N}}
\newcommand{\Ss}{\mathcal{S}}
\newcommand{\VV}{\mathcal{V}}

\newcommand{\KK}{\mathcal{K}}
\newcommand{\CC}{\mathcal{C}}
\newcommand{\LL}{\mathcal{L}}
\newcommand{\TT}{\mathcal{T}}
\newtheorem{lemma}{Lemma}

\newlength\fheight
\newlength\fwidth

\title{Exact and Heuristic Methods with Warm-start for Embedded Mixed-Integer Quadratic Programming Based on Accelerated Dual Gradient Projection
	  }

\author{
	Vihangkumar V. Naik\thanks{Vihangkumar V. Naik is with ODYS S.r.l., Via A. Passaglia 185, Lucca, Italy. Email: \texttt{vihang.naik@odys.it}.} \and 
	Alberto Bemporad\thanks{A. Bemporad is with the IMT School for Advanced Studies Lucca, Italy. Email: \texttt{alberto.bemporad@imtlucca.it}.}
}

\begin{document}
	\maketitle
	
\begin{abstract}
Small-scale Mixed-Integer Quadratic Programming (MIQP) problems often
arise in embedded control and estimation applications. Driven by the need for 
algorithmic simplicity to 
target computing platforms with limited memory and computing resources, 
this paper proposes a few approaches to solving MIQPs, either to optimality or
suboptimally. We specialize an existing Accelerated Dual Gradient Projection (GPAD) 
algorithm to effectively solve the Quadratic Programming (QP) relaxation that arise during Branch and Bound (B\&B) and propose a generic framework to warm-start the binary variables which reduces the number of QP relaxations.
Moreover, in order to find an 
integer feasible combination of the binary variables upfront, two heuristic approaches are presented: 
($i$) without using B\&B, and ($ii$) using B\&B with a significantly reduced number of QP relaxations.
Both heuristic approaches return an integer feasible solution that may be suboptimal but involve a much reduced computation effort. Such a feasible solution can be either implemented directly or used to set an initial upper bound on the optimal cost in B\&B.
Through different hybrid control and estimation examples involving binary decision variables, we show that the performance of the proposed methods, although very simple to code, is comparable to that of state-of-the-art MIQP solvers.
\end{abstract}


\section{Introduction}

We consider a Mixed-Integer Quadratic Programming (MIQP) problem of the following general form
\begin{subequations}\label{eq:QP}
	\begin{align}
	\displaystyle{\min_{z}}\quad
	& V(z)\eqdef \frac{1}{2}z'Qz+c'z\label{eq:QP1}\\ 
	\st\quad            & \ell\leq Az\leq u            \label{eq:QP2}\\
	& A_{eq}z=b_{eq} \label{eq:QP3}\\
	& \bar A_iz\in\{\bar\ell_i,\bar u_i\},\ i=1,\ldots,p\label{eq:QP4}
	\end{align}	
\end{subequations}
where $z\in \mathbb{R}^n$ is the vector of decision variables, $Q=Q'\succ 0\in \mathbb{R}^{n \times n}$ is the Hessian matrix, 
$c\in \mathbb{R}^n$, $A\in \mathbb{R}^{m \times n}$, $\ell,u\in \mathbb{R}^m$, $\ell \leq u$, represents linear inequality constraints, $A_{eq} \in \mathbb{R}^{q \times n}$, $b_{eq} \in \mathbb{R}^q$ describes linear equality constraints and the integrality constraints are described by $\bar{A} \in \mathbb{R}^{p \times n}$, $\bar{\ell},\bar{u} \in \mathbb{R}^{p}$, 
$\bar\ell \leq \bar u$, and $p\leq m$.

More frequently, MIQP problems are expressed by just restricting some of the optimization variables to be binary. This is a special case of~\eqref{eq:QP4}, as binary constraints $z_i\in\{0,1\}$
for $i=1,\ldots,p$, $p\leq n$, correspond to setting $\bar \ell_i=0, \bar u_i=1$, and $\bar A_i$ as
the $i$-th row of the identity matrix. 

MIQP problems arise in many applications such as 
hybrid Model Predictive Control (MPC)~\cite{BM99}, moving-horizon estimation~\cite{BMM99c,
FMM02}, piecewise affine regression~\cite{Bemporad_roll_hyb_id_cdc2001,
roll2004identification, MNPB18SYSID}, trajectory generation~\cite{mellinger2012mixed},
economic dispatch~\cite{papageorgiou2007mixed},
planning and design~\cite{propato2004booster,yang2007tcsc}, and
scheduling~\cite{catalao2010scheduling}.
In particular, MPC based on hybrid dynamical models has been adopted in various fields, due to its capability to handle process dynamics coupled with logical rules, switching dynamics, discrete actuation signals, and mixed linear and logical constraints. 
On-line implementation of hybrid MPC requires the solution of an MIQP problem of the form~\eqref{eq:QP} at every sampling instant~\cite{BM99}. 
The Binary Quadratic Programming (BQP), which is a special case of~\eqref{eq:QP} in which all the decision variables are binary, was used for example in~\cite{MNPB18CDC_energy} to solve an energy disaggregation problem catering to the real-time implementation on a smart meter. 

The computational complexity of MIQP in general grows exponentially with the number of integrality constraints. 
This inherent characteristic usually restricts the solver to run on a desktop computer, for which excellent commercial software packages exist~\cite{Cplex14,gurobi14,XPRESS-FICO,mosek}. 
However, due to their complexity and unavailability of library-free source code, 
these packages are not amenable for implementation in embedded platforms. Hence,
how to solve small-scale MIQP problems by means of simpler algorithms has attracted considerable attention from the scientific community in recent years. 

Various approaches have been proposed to solve MIQP problems to optimality, mostly based on the well-known branch and bound (B\&B) algorithm~\cite{floudas1995nonlinear}. B\&B relies on solving a sequence of Quadratic Programming (QP) relaxations
\begin{subequations}
    \begin{align}
    \displaystyle{\min_{z}}\quad
    & V(z)\eqdef \frac{1}{2}z'Qz+c'z\label{eq:QP1_bnb}\\ 
    \st\quad            & \ell\leq Az\leq u            \label{eq:QP2_bnb}\\
    & A_{eq}z=b_{eq} \label{eq:QP3_bnb}\\
    &\label{eq:AIuz=bIu}\bar{A}_{\II_{\bar u}}z=\bar{u}_{\II_{\bar u}}\\
    &\label{eq:AIellz=bIell}\bar{A}_{\II_{\bar \ell}}z=\bar{\ell}_{\II_{\bar \ell}}\\ 
    &\label{eq:ellJ<=AJz} \bar{\ell}_{\JJ} \leq \bar{A}_{\JJ}z \leq \bar{u}_{\JJ}
    \end{align} \label{eq:QP_bnb}
    \label{eq:QP_relaxation}
\end{subequations}
where $\II_{\bar u},\II_{\bar \ell}\subseteq\{1,\ldots,p\}$,
$\II_{\bar u}\cap \II_{\bar \ell} = \emptyset$, 
$\JJ=\{1,\ldots,p\} \setminus (\II_{\bar u}\cup \II_{\bar \ell})$.
The authors in~\cite{axehill_and_hansson_2006mixed} proposed an approach using a dual active-set method for QP within B\&B for online solution MIQPs arising in hybrid model predictive control (MPC). However, this approach does not exploit dual lower bounds on the optimal cost, that is very useful to terminate the QP solver prematurely~\cite{fletcher_and_Leyffer1998numerical, axehill2008dual}. 
Recently, various B\&B methods for solving MIQPs in embedded platforms have been presented, which mainly differ
for how QP relaxations~\eqref{eq:QP_relaxation} are solved, such as interior-point methods~\cite{frick2015embedded}, active-set method based on nonnegative least squares (NNLS)~\cite{bemporad_miqpnnls2015}, an active-set method based on NNLS and proximal-point iterations~\cite{BN18NMPC_nnlsprox_miqp}, and Alternating Direction Method of Multipliers (ADMM) based on the Operator Splitting Quadratic Program (OSQP) solver~\cite{stellato2018}.
The latter two methods have an advantage that they do not require regularization of the Hessian matrix, which might instead introduce a bias on the optimal solution. For small problems having few binary variables, multi-parametric programming~\cite{bemporad2002explicit,Bem15,AB09,BBBM05} provides
an alternative to embed the solver code in the application, as~\eqref{eq:QP} 
is pre-solved offline to get the optimizer $z^*$ as a function
of parameters possibly entering $c$, $\ell$, $u$, $b_{eq}$ in a linear fashion.

The aforementioned MIQP solution methods were tailored to solving small-scale MIQPs, such as those typically arising in embedded applications. In an embedded environment, several restrictions exist on available memory, throughput, worst-case execution time, and software simplicity for code analysis, that translate
into restrictions on the MIQP algorithm. 
A number of embedded QP solvers have been developed, including the active-set method qpOASES~\cite{Ferreau2014}, NNLS~\cite{bemporad_ieeetac-qpnnls, Bem18}, the interior-point method CVXGEN~\cite{Mattingley2012},  FORCES~\cite{Domahidi_forces2012}, the first-order method FiOrdOs~\cite{firordos}, accelerated dual gradient projection (GPAD)~\cite{patrinosBemporad2014GPAD},  OSQP~\cite{osqp}, and commercial implementations for industrial production are even available, such as ODYS QP solver~\cite{ODYSQP}.

Among QP solution methods, Nesterov's fast gradient method~\cite{nesterov1983method, nesterov2003introductory} has received great attention recently within the embedded control community for embedded applications~\cite{Richter_fgm_mpc_2011, kogel2011fast, patrinosBemporad2014GPAD, Giselsson_fgm_mpc_2014}. The primary reasons for this are its simplicity, relatively good performance, and good bounds on the worst-case number of iterations. Such characteristics make the method favorable for real-time embedded applications. The fast gradient method has been employed for a field-programmable gate array (FPGA) implementation of embedded MPC~\cite{jerez2013embedded} and for using reduced-precision fixed-point arithmetic~\cite{jerez_megahertz_2014}.
An accelerated gradient projection method applied to solve the dual QP problem (GPAD) was proposed by~\cite{patrinosBemporad2014GPAD}. For the dual gradient projection algorithm, the fixed-point implementation~\cite{patrinos2013fixed} and FPGA implementation have been studied in~\cite{rubagotti2016real}. 

In spite of the abundance of efficient algorithms available to solve QPs as in~\eqref{eq:QP_relaxation}, solving an MIQP to optimality may be impractical in fast applications due to the integrality constraints~\eqref{eq:QP4}. For this reason, efforts have been made recently for solving MIQPs \emph{approximately} in a quick manner. Heuristic methods were proposed which lead to a suboptimal solution of the MIQP while taking considerably less computation time. 
A heuristic method for finding a feasible solution of a generic mixed-integer program, known as the \emph{feasibility pump}, was proposed in~\cite{fischetti2005feasibility}. Further,~\cite{bertacco2007feasibility} extended this approach to binary and general-integer variables, see also~\cite{achterberg2007improving}~and~\cite{heu_MIP_FischettiLodi}.
Recently,~\cite{TMBB17} presented a heuristic based on quantization to find an approximate solution of MIQP problem via ADMM.
An approach based on an operator-splitting method for finding suboptimal solutions, with guaranteed local convergence results under certain assumptions was presented in~\cite{frick_lowcomplexityhmpc}.

By largely extending the original idea presented in~\cite{NBIFAC17},
this paper proposes an embedded MIQP solution method based on a 
specialized version of the GPAD algorithm from~\cite{patrinosBemporad2014GPAD} 
to solve the relaxed QP subproblems arising during B\&B. 
The GPAD algorithm is equipped with an infeasibility detection mechanism, which is essential to detect QP subproblems with infeasible binary combinations that occurs frequently during B\&B.
As during B\&B, each QP subproblem changes from the parent one by only turning an inequality constraint to an equality one, 
using a dual QP method like~\cite{patrinosBemporad2014GPAD} has the advantage that this simply corresponds to removing the sign restriction on the corresponding dual variable. Therefore, the preconditioning required by the gradient projection algorithm
can be executed just once at the root node and applied to the entire set of QP relaxations generated during B\&B. In addition, the optimal dual solution of the parent problem and current branching variable information are used to warm-start the two following children QP subproblems. 
Moreover, while solving relaxed QP subproblems, since each iterate of the GPAD algorithm is dual feasible due to the projection operation, it is possible to terminate
prematurely the QP solver in case it leads to a cost that is greater than the current best known upper-bound on the MIQP solution and, accordingly, prune the B\&B tree.
Besides preconditioning, we also use restart~\cite{o2015adaptive} to enhance the overall performance of the dual QP solver.

In this paper we also provide a generic framework to warm-start binary variables, following
preliminary ideas recently introduced by the authors in the conference paper~\cite{BN18NMPC_nnlsprox_miqp} (such ideas have been already adopted in~\cite{Hespanhol2019_struc_exploit_bnb, marcucci2020warm}). We further show that, by using a different branching strategy, one may reduce the number of solved QP relaxations. The warm-starting framework is also extended for efficient solution of MIQP problems where the binary variables are subject to an exclusive-or condition, such as due to one-hot encoding of categorical decision variables. For MIQPs generated by receding-horizon optimization (such as hybrid MPC, moving horizon estimation, and piecewise affine regression),
the framework allows warm starts generated by the optimal values calculated at the previous sampling instant or by artificial neural networks~\cite{MB2019_learning_ws}. 
Furthermore, this paper presents two heuristic approaches tailored to providing a good upper bound on the optimal cost for initializing the B\&B algorithm.
The first heuristic is a simple modification of the GPAD algorithm and requires no branch and bound. 
The second method reduces the number of binary variables to branch upon during B\&B by fixing the remaining binaries to one of the bounds. 
Though there is no guarantee of convergence for the heuristic methods, in practice they turn out to be quite effective in solving MIQPs approximately, in most cases very close to the optimal solution. 
The heuristic approaches are also combined with the warm-start framework for binary variable as a special case. 
The effectiveness of the proposed approaches are demonstrated in numerical examples.

\subsection{Notation}

Let $\rr^n$, $\rr^{m \times n}$, $\N$ denote the set of real vectors of length $n$, the set of real matrices of dimension $m$ by $n$, and
the set of natural integers, respectively. Let $\II\subset\N$ be a finite set of integers
and denote by $\card(\II)$ its cardinality. 
For a vector $a\in\rr^{n}$, $a_i$ denotes the $i$-th entry of $a$, 
$a_\II$ the subvector obtained by collecting the entries $a_i$ for all $i\in\II$, 
$\|a\|_\infty$ the infinite norm of $a$, 
the condition $a>0$ is equivalent to 
$a_i>0$, $\forall i=1,\ldots,n$ (and similarly for $\geq$, $\leq$, $<$).
We denote by $0_n$ the vector of $\rr^n$ with all zero components, with the subscript $_n$
dropped whenever the dimension is clear from the context.
For a matrix $A\in\rr^{m\times n}$, its transpose is denoted by $A'$, its
$i$th row by $A_i$, 
the submatrix of $A$ obtained by collecting $A_i$ for all $i\in\II$ is denoted by $A_\II$, and its Frobenius norm by $\|A\|_F$, its maximum eigenvalue by  $\lambda_{\rm max}(A)$,
and its maximum singular value by $\sigma_{\rm max}(A)$.
Matrix $A^{-1}$ denotes the inverse of a square matrix $A\in\rr^{n\times n}$ (if it exists),
$A\succ0$ denotes that $A$ is positive definite, 
and similarly $\succeq$, $\prec$, $\preceq$ denote positive semidefiniteness, 
negative definiteness, negative semidefiniteness, respectively. 
$A_{ij}$ represent the element $(i,j)$ of matrix $A$.

\section{Accelerated dual gradient projection (GPAD)}
\label{sec:GPAD}

The dual problem of the QP~\eqref{eq:QP_relaxation}, where for simplicity of notation
we consider the case $\II_{\bar u}=\II_{\bar \ell}=\emptyset$, is the following convex QP 
\begin{subequations}\label{eq:dual-QP}
	\begin{align}
	\displaystyle{\max_{\lambda,\nu}}~
	& \Psi(\lambda,\nu)\eqdef-\frac{1}{2}\smallmat{\lambda\\\nu}'
	\Aa Q^{-1}\Aa'\smallmat{\lambda\\\nu}-\DD'\smallmat{\lambda\\\nu}-\frac{1}{2}c'Q^{-1}c\nonumber\\
	\st\quad & \lambda\geq 0,\ \nu~\mbox{free} \label{eq:dual}
	\end{align}
	where 
	\begin{equation}
\Aa=\smallmat{A\\-A\\\bar A\\-\bar A\\A_{eq}},\
\BB=\smallmat{u\\-\ell\\\bar u\\-\bar \ell\\b_{\rm eq}},\
\DD=\smallmat{d_u\\d_\ell\\\bar d_{ u}\\\bar d_{\ell}\\f}=\BB+\Aa Q^{-1}c
	\label{eq:df}
\end{equation}
$\lambda=[\lambda_u'\ \lambda_\ell'\ \bar \lambda_{u}'\ \bar \lambda_{\ell}']'$
is the vector of dual variables associated with inequality constraints,
	$\lambda_{u},\lambda_{\ell},d_{u},d_{\ell}\in\rr^m$,
	$\bar \lambda_{u},\bar \lambda_{\ell},\bar d_{ u},\bar d_{\ell}\in\rr^p$,
$\nu$ is the vector of dual variables corresponding to equality constraints,
and $\nu,f\in\rr^q$.
\end{subequations}

The accelerated dual gradient projection (GPAD) method described in
Algorithm~\ref{algo:GPAD} was proposed by~\cite{patrinosBemporad2014GPAD} to solve a strictly convex QP problem as in~\eqref{eq:QP_relaxation} by applying Nesterov's fast gradient method to its dual problem. 

\begin{algorithm}[tb!]
	\caption{Accelerated gradient projection method to solve dual QP problem~\eqref{eq:dual-QP}
}\label{algo:GPAD}
~~\textbf{Input}: matrices $Q$, $A$, $A_{eq}$, $\bar A$, and vectors $c$, $\ell$, $u$, $b_{eq}$,
	$\bar\ell$, $\bar u$.
	\vspace*{.1cm}\hrule\vspace*{.1cm}
	
	\begin{enumerate}[label=\arabic*., ref=\theenumi{}]
	\item $H \leftarrow \Aa Q^{-1}\Aa'$, $L\leftarrow\|H\|_F$
		or $L\leftarrow\lambda_{\rm max}(H)$,\newline
		$\Aa_L\leftarrow \frac{1}{L}\Aa$, $\BB_L\leftarrow \frac{1}{L}\BB$; \label{step:GPAD init}
\item $\lambda_{0},\lambda_{-1}\leftarrow0$, $\nu\leftarrow0$;		 
	\item $k\leftarrow0$; \label{step:k0}
		\item \textbf{repeat} \label{step:repeat}
		\begin{enumerate}[label=\theenumi{}.\arabic*., ref=\theenumi{}.\theenumii{}]
		\item $\beta_{k}\leftarrow\max\left\{\frac{k-1}{k+2},0\right\};\label{step:beta}$
		\item $\smallmat{w_{k}\\w_{eq,k}} \leftarrow  \smallmat{
			\lambda_{k} \\ \nu_{k}} + \beta_{k} 
		\left(\smallmat{
			\lambda_{k} \\  \nu_{k}
		}  - 
		\smallmat{\lambda_{k-1} \\  \nu_{k-1}} \right); $\label{eq:gpad steps w_k} 
		\item $z_{k} \leftarrow-Q^{-1}\Aa'\smallmat{w_{k}\\w_{eq,k}}-Q^{-1}c;$ \label{eq:gpad steps z_k}
		\item $\smallmat{s_{k}\\s_{eq,k}} \leftarrow 
		\left (\Aa_L z_{k}-\BB_L
		\right );$ \label{eq:gpad steps s_k}
		\item ${\lambda_{k+1}}\leftarrow \mathrm{max}\{{w_{k}}+{s_{k}},0\}$;
		\label{step:yk1}	
		\item $\nu_{k+1}\leftarrow w_{eq,k}+s_{eq,k}$;\label{step:nuk1}
        \end{enumerate}
		\item [] \textbf{until} convergence;
		\item $z^*\leftarrow z_k$, $\lambda^*\leftarrow w_k$, $\nu^* \leftarrow w_{eq,k}$;
		\item $a^*\leftarrow \Aa'\smallmat{\lambda^*\\\nu^*}$, $V^*\leftarrow-\frac{1}{2}(a^*)'Q^{-1}a^*-\BB'\smallmat{\lambda^*\\\nu^*}
		-(Q^{-1}c)'(a^*+\frac{1}{2}c)= \Psi^*.$ 
    	\end{enumerate}
	
	\vspace*{.1cm}\hrule\vspace*{.1cm}
	~~\textbf{Output}: Primal solution $z^*$, 
		optimal cost $V^*$, dual solution $(\lambda^*,\nu^*)$.
\end{algorithm}

Note that the only difference between inequality and equality constraints is
that Step~\ref{step:yk1} of Algorithm~\ref{algo:GPAD} projects the dual variables 
corresponding to inequalities on the nonnegative orthant, while Step~\ref{step:nuk1} does not for equality constraints. We will exploit this in the B\&B approach described in Section~\ref{sec:BB} and in the heuristic approaches presented in Section~\ref{sec:heuristics}.

The Lipschitz constant $L$ of the gradient of the dual QP function $\Psi$ plays an important role in Algorithm~\ref{algo:GPAD}, as it affects its convergence speed~\cite{patrinosBemporad2014GPAD}. A valid value of $L$
is the Frobenius norm $\|H\|_F$ of the dual Hessian $H=\Aa Q^{-1}\Aa'$. In case a Cholesky factorization $Q=\CC'\CC$ is used to compute $Q^{-1}$, a better (smaller or equal) value for $L$ is $\lambda_{\rm max}(H)$, which can be evaluated by computing the largest
singular value $\sigma_{\rm max}$ of matrix $\Aa\CC^{-1}$,
as $\sigma^2_{\rm max}(\Aa\CC^{-1})=\lambda_{\rm max}(\Aa\CC^{-1}(\Aa\CC^{-1})')=
\lambda_{\rm max}(\Aa Q^{-1}\Aa')$.

\subsection{Stopping criteria}

The iterations of Algorithm~\ref{algo:GPAD} are terminated when the primal feasibility criterion is satisfied along with an optimality criterion~\cite{patrinosBemporad2014GPAD}. 
The primal feasibility criterion is
\begin{align}
&s_{k}^{j} \leq \frac{1}{L}\epsilon_{G},~\forall j = 1, \ldots, 2(m+p) \nonumber \\
&|s_{eq,k}^{j}| \leq \frac{1}{L}\epsilon_{G},~\forall j = 1, \ldots, q 	\label{eq:gpad feas cond}
\end{align}
where the feasibility tolerance is $\epsilon_{G}>0$. 
A sufficient condition that can be used as an optimality stopping criterion is
\begin{equation}
 -\smallmat{w_{k}\\w_{eq,k}}'\smallmat{s_{k}\\s_{eq,k}} \leq \frac{1}{L}\epsilon_{V}, \quad 	w_{k}\geq0
\label{eq:gpad opt cond}
\end{equation}
where the optimality tolerance is $\epsilon_{V} \geq0$, which is derived from the duality gap calculation  $V(z_{k})-V^{*} \leq V(z_{k})-\Psi(\smallmat{w_{k}\\w_{eq,k}}) = -\smallmat{w_{k}\\w_{eq,k}}'\smallmat{s_{k}\\s_{eq,k}} \leq \frac{1}{L}\epsilon_{V}$.

\subsection{Preconditioning}

It is well known that the convergence speed to the solution of first-order optimization methods can be often largely improved by \emph{preconditioning} the matrices defining the problem~\cite{Ber99, giselsson2014preconditioning}. We use the Jacobi diagonal scaling~\cite{Ber09} applied to the dual Hessian matrix $H=\Aa Q^{-1}\Aa'$,
defined by replacing
\begin{subequations}
	\begin{align}
	\Aa_j &\leftarrow \bar \theta_{j}\Aa_j,\ \BB_j\leftarrow \bar \theta_j\BB_j\\
	\bar \theta_j&\eqdef\frac{1}{\sqrt{\Aa_j Q^{-1}\Aa_j'}}\\
    &j=1,\ldots,2(m+p)+q\nonumber
	\end{align}
	\label{eq:scaling}
\end{subequations}

\subsection{Restart}

Accelerated gradient projection methods often exhibit non-monotonic decrease of the objective
function, with ripples observed in the sequence of objective values. Employing a scheme to restart
the sequence of $\beta_k$ in Algorithm~\ref{algo:GPAD} 
can largely improve the convergence property of the method. We use the gradient-based adaptive restart idea of~\cite{o2015adaptive,giselssonBoyd_restart} for the dual problem~\eqref{eq:dual-QP} by checking the following condition
\begin{equation}
-\nabla \Psi\left(\smallmat{w_{k}\\w_{eq,k}}\right)'\left(\smallmat{\lambda_{k+1}\\\nu_{k+1}}-\smallmat{\lambda_{k}\\\nu_{k}}\right)>0 \label{eq:restart condition}
\end{equation}
Whenever condition~\eqref{eq:restart condition} is satisfied, the value of scalar $k$ in Step~\ref{step:beta} is zeroed, which in turn resets the momentum term $\beta_k$ to zero. 
The main advantage of the gradient-based restart condition~\eqref{eq:restart condition} is that it can be immediately computed by available quantities.

\subsection{Infeasibility detection}    

Infeasibility detection in first-order methods has been studied in~\cite{raghunathanCairano_Infeasibility, O'donoghue:2016:COV:2944533.2944554,osqp-infeasibility}.     
During B\&B, equality constraints are generated be imposing that
$\bar A_iz$ is either $\bar \ell_i$ or $\bar u_i$ for a subset of the indices.
The resulting QP problem~\eqref{eq:QP_relaxation}
may be infeasible, which makes the dual cost $\Psi(\lambda_k,\nu_k)$ diverge to $+\infty$. 
The following Lemma~\ref{lemma:infeas} characterizes the asymptotic
behavior of Algorithm~\ref{algo:GPAD} in case of infeasibility of~\eqref{eq:QP_relaxation}.

\begin{lemma}
\label{lemma:infeas}
Let the QP problem~\eqref{eq:QP_relaxation}
be infeasible.
Then $\lim_{k\rightarrow\infty}\Aa'\smallmat{\lambda_k\\\nu_{k}}/
\|\smallmat{\lambda_k\\\nu_{k}}\|_\infty=0$.
\end{lemma}

\proof 
Let $\eta_k=\smallmat{\lambda_k\\\nu_{k}}$. Since $H=\Aa Q^{-1}\Aa'\succeq 0$, 
for $\Psi(\lambda_k,\nu_{k})\rightarrow+\infty$ it must occur that $\DD'\eta_k\rightarrow -\infty$, and therefore some components $\eta_{k,i}$ must tend to $+\infty$ for corresponding negative entries in vector $\DD$. 
Assume now by contradiction that the quantity $\Aa'\frac{\eta_k}{\|\eta_k\|_\infty}$ does not go to zero for $k\rightarrow\infty$. In this case
there would exist a subsequence $\eta_r$ such that $
\|\Aa'\frac{\eta_r}{\|\eta_r\|_\infty}\|_2\geq\epsilon$ for some $\epsilon>0$.
In this case, $\Psi(\lambda_r,\nu_r)=-\frac{1}{2}\eta_r'\Aa Q^{-1}\Aa'\eta_r-\DD'\eta_r\leq-\frac{1}{2\lambda_{\rm max}(Q)}\|\Aa'\eta_r\|_2^2+\|\DD\|_2\|\eta_r\|_2\leq
-\frac{\epsilon^2}{2\lambda_{\rm max}(Q)}\|\eta_r\|_\infty^2+\|\DD\|_2\sqrt{2(m+p)+q}\|\eta_r\|_\infty$, where $\lambda_{\rm max}(Q)$ is the largest eigenvalue of $Q$. Since some components of $\eta_r$ diverge, $\|\eta_r\|_\infty$
diverges as well, and therefore $\Psi(\lambda_r,\nu_{r})\leq 0$ for $r$ sufficiently large. This contradicts the fact that $\Psi(\lambda_k,\nu_{k})\rightarrow+\infty$,
and therefore $\Aa'\frac{\eta_k}{\|\eta_k\|_\infty}$ must go to zero
asymptotically if the QP is infeasible.
\hfill\QED

Motivated by Lemma~\ref{lemma:infeas}, we propose the infeasibility detection criterion summarized in Algorithm~\ref{algo:Infeasibility detection},
where $\epsilon_{I}> 0$ is a given infeasibility detection tolerance.

\begin{algorithm}[ht!]
	\caption{Infeasibility detection}\label{algo:Infeasibility detection}
	\begin{enumerate}[label=\arabic*., ref=\theenumi{}]
		\item $\alpha_k\leftarrow\|\smallmat{\lambda_k\\\nu_{k}}\|_\infty$; 
		\item \textbf{if} {$\| \Aa'\smallmat{\lambda_k\\\nu_{k}}\|_\infty \leq \epsilon_I \alpha_k$ \textbf{and} $\DD'\smallmat{\lambda_k\\\nu_{k}} < -\epsilon_I \alpha_k$}\label{step:infeas}
		\textbf{then stop} (problem is infeasible)
	\end{enumerate}
\end{algorithm}
 
By letting $\eta_k=\smallmat{\lambda_k\\\nu_{k}}$, $\mu_k=\frac{\lambda_k}{\|\eta_k\|_\infty}$, $\pi_k=\frac{\nu_{k}}{\|\eta_k\|_\infty}$,
$\mu_k\in\rr^{2(m+p)}$, $\pi_k\in\rr^q$, the criterion in Step~\ref{step:infeas} of Algorithm~\ref{algo:Infeasibility detection} amounts to verify the following conditions
\begin{align}
	& \left\{\begin{matrix}
	&\Aa'\smallmat{\mu_k\\\pi_k}\approx 0\\[.2em]
	&\BB'\smallmat{\mu_k\\\pi_k}<0\\[.2em]
	&\mu_k\geq 0
	\end{matrix}\right.
	\label{eq:Farkas}
\end{align} 
According to Farkas Lemma~\cite[p.~201]{Roc70}, condition~\eqref{eq:Farkas} is equivalent to an indication that the QP~\eqref{eq:QP_relaxation} is infeasible.
Moreover, $b'\smallmat{\mu_k\\\pi_k}<0$ can be equivalently replaced by $d'\smallmat{\mu_k\\\pi_k}<0$ when $\Aa'\smallmat{\mu_k\\\pi_k}=0$.

\subsection{Early stopping criterion for the objective function}

As we will describe in Section~\ref{sec:BB}, when an upper bound $V_0$ on the MIQP optimal cost is available from an existing integer feasible solution, a node of the B\&B tree can be fathomed when the QP relaxation 
is discovered to have a cost greater than or equal to $V_0$. 
We therefore immediately stop Algorithm~\ref{algo:GPAD} whenever 
\begin{equation}
\Psi(\lambda_k,\nu_k) \geq V_{0}
\label{eq:V0}
\end{equation}
as $V(z_k)\geq \Psi(\lambda_k,\nu_k)\geq V_0$.
Since the quantity $a_k=\Aa'\smallmat{\lambda_k\\\nu_k}$ is already available from Step~\ref{step:infeas} of the infeasibility detection Algorithm~\ref{algo:Infeasibility detection}, the quantity $\Psi(\lambda_k,\nu_k)$ can be easily computed as 
\begin{align}
&\Psi(\lambda_k,\nu_k)=-\frac{1}{2}a_k'Q^{-1}a_k-\BB'\smallmat{\lambda_k\\\nu_k}-(Q^{-1}c)'(a_k+\frac{1}{2}c)
\end{align}

In order to best use the early stopping criterion~\eqref{eq:V0}, in Section~\ref{sec:heuristics} we propose heuristic approaches to compute an integer feasible solution providing the required upper bound $V_0$ upfront.

\section{Branch and Bound MIQP Algorithm}
\label{sec:BB}

Algorithm~\ref{algo:MIQP-GPAD} proposes a B\&B scheme for solving the MIQP problem~\eqref{eq:QP}
based on Algorithm~\ref{algo:GPAD} for solving QP relaxations. We start with formalizing the standard B\&B algorithm followed by our specific contributions. As it is common, the B\&B procedure
can be described by a search tree where each node represents a unique QP subproblem. 
The tree is initialized by relaxing all the integrality constraints \eqref{eq:QP4}, allowing them to take any values in the interval $\left[\bar{\ell}, \bar u \right] $. The 
corresponding relaxed QP problem~\eqref{eq:QP_relaxation} with $\II_{\bar u}=\II_{\bar \ell}=\emptyset$
is termed as the \emph{root} node of the tree. 
During \emph{branching}, two \emph{children} nodes are created from each \emph{parent} node
by fixing $\bar{A}_{i} z = \bar{\ell}_i$ and $\bar{A}_{i} z = \bar{u}_i$, respectively. Nodes with all 
integrality constraints fixed to either $\bar{\ell}_i$ or $\bar{u}_i$
($\II_{\bar u}\cup\II_{\bar \ell}=\{1,\ldots,p\}$) are called \emph{leaf} nodes of the tree. 
Let $\TT$ be a unique tuple containing the indices $\II_{\bar \ell}$, $\II_{\bar u}$ associated with a QP
relaxation. The set $\Ss$ is a stack which stores the tuples $\TT$ of pending relaxations
to be solved by Algorithm~\ref{algo:GPAD}.

At Step~\ref{step:pop}, the last tuple $\TT$ is popped from the stack $\Ss$ and the corresponding QP relaxation~\eqref{eq:QP_bnb} is solved at Step~\ref{step:QP} using Algorithm~\ref{algo:GPAD}. The additional stopping criterion~\eqref{eq:V0} is included in Algorithm~\ref{algo:GPAD}, where $V_0$ denotes cost of the best integer-feasible solution found so far.

Step~\ref{step:QPsolved} checks if the solution of the QP relaxation is feasible and if the optimal solution $V^*$ of the QP subproblem satisfies the condition $V^*\leq V_0$. 
If the condition at Step~\ref{step:QPsolved} is satisfied then Step~\ref{step:check-integrality} checks
if all integrality constraints in~\eqref{eq:QP4} are satisfied and subsequently updates the best known integer-feasible cost $V_0$ and the corresponding vector of decision variables $\zeta^*$. Otherwise, branching  is carried out by executing Steps~\ref{step:branch-start}--\ref{step:branch-end}. At Step~\ref{step:branch-start}, a branching variable $j$ is selected from the set $\JJ$ such that $t_i=\bar A_iz^*$ is most distant from $\bar\ell_i$ and $\bar u_i$. The index $j$ is moved from the set of inequality constraints $\JJ$ into sets of equality constraints and corresponding two new QP subproblems $\TT_0$, $\TT_1$ are formed at Step~\ref{step:AA}.
These two problems are pushed onto the stack $\Ss$, giving higher priority to the problem with the least fractional part at Step~\ref{step:branch-end}.

Once the stack $\Ss$ becomes empty, no further QP relaxations must be solved. At Step~\ref{step:Vstar}, if the value of $V_0$ is still $+\infty$ then the given MIQP problem~\eqref{eq:QP} is declared infeasible; otherwise, the optimal solution $\zeta^*$ and optimal cost $\VV^*$ are returned. 

\begin{algorithm}[tb!]
	\caption{MIQP solver based on GPAD}
	\label{algo:MIQP-GPAD}
	~~\textbf{Input}: MIQP problem matrices $Q=Q'\succ0$, $A$, $\bar A$, $A_{eq}$ and vectors c, $\ell$, $u$, $b_{eq}$, $\bar\ell$, $\bar u$; tolerance values $\epsilon_{G} > 0$, $\epsilon_{V}\geq 0$, $\epsilon_{I}> 0$.
	\vspace*{.1cm}\hrule\vspace*{.1cm}
	
	\begin{enumerate}[label=\arabic*., ref=\theenumi{}]
		
		\item \label{step:init}\textbf{set} 
		$V_0\leftarrow +\infty$; $\zeta^*\leftarrow \emptyset$;
		$\II_{\bar \ell}\leftarrow \emptyset;\II_{\bar u}\leftarrow \emptyset$; 
		$\TT\leftarrow(\II_{\bar \ell},\II_{\bar u})$;
		$\Ss\leftarrow\{\TT\}$;
		\item \textbf{while} $\Ss\neq\emptyset$ \textbf{do}:
		
		\begin{enumerate}[label=\theenumi{}.\arabic*., ref=\theenumi{}.\theenumii{}]
			\item \label{step:pop} \textbf{pop} $\TT=(\II_{\bar \ell},\II_{\bar u})$ from $\Ss$  ($\Ss\leftarrow\Ss\setminus\{\TT\}$);
            \item $\JJ\leftarrow\{1,\ldots,p\}\setminus(\II_{\bar \ell}\cup \II_{\bar u})$;
			\item \label{step:QP}\textbf{execute} Algorithm~\ref{algo:GPAD}
			 to solve~\eqref{eq:QP_bnb}
			 under condition~\eqref{eq:V0};
			\item \label{step:QPsolved}\textbf{if} the solution $z^*,V^*$ is returned \textbf{and} $V^* \leq V_0$ \textbf{then}
			
			\begin{enumerate}[label=\theenumii{}.\arabic*., ref=\theenumii{}.\theenumiii{}]
				\item \label{step:check-integrality}\textbf{if} $\JJ=\emptyset$ 
                \textbf{or} 
				$t_i\eqdef\bar A_iz^*\in\{\bar \ell_i,\bar u_i\}$, 
				$\forall i\in \JJ$ 
				\textbf{then}
				$V_0\leftarrow V^*$, $\zeta^*\leftarrow z^*$; \textbf{otherwise}
				
				\begin{enumerate}[label=\theenumi{}.\theenumii{}.\arabic*., ref=\theenumi{}.\theenumii{}.\theenumiv{}]
					\item \label{step:branch-start}$\displaystyle{j\leftarrow\arg\min_{i\in
							J}}$ $\left|t_i-\frac{\bar \ell_i+ \bar u_i}{2}\right|$; 
					\item \label{step:AA} 
                    $\TT_0\leftarrow(\II_{\bar \ell}\cup\{j\},\II_{\bar u})$;
					$\TT_1\leftarrow(\II_{\bar \ell},\II_{\bar u}\cup\{j\})$;
					\item \label{step:branch-end}\textbf{if} $t_{j}\leq\frac{\bar \ell_i+\bar u_i}{2}$ \textbf{then} push $\TT_1$ and then $\TT_0$ on $\Ss$ \textbf{otherwise} push $\TT_0$ and then $\TT_1$;
				\end{enumerate}
			\end{enumerate}
		\end{enumerate}
		\item \label{step:Vstar}\textbf{if} $V_0=+\infty$ \textbf{then}~\eqref{eq:QP} infeasible
		\textbf{otherwise}  
		$\VV^*\leftarrow V_0$;
		\item \textbf{end}.
	\end{enumerate}
	
	\vspace*{.1cm}\hrule\vspace*{.1cm}
	~~\textbf{Output}: Solution $\zeta^*$ of the MIQP problem~\eqref{eq:QP}, optimal cost $\VV^*$, or infeasibility status.
\end{algorithm}

\subsection{Exploiting the fixed structure of dual QP relaxations}

During the execution of the B\&B algorithm,
from one QP subproblem to another only the constraints~\eqref{eq:AIuz=bIu}-\eqref{eq:ellJ<=AJz} are changed. 
Such changes simply map into the following
conditions imposed at Steps~\ref{step:yk1}--\ref{step:nuk1} of Algorithm~\ref{algo:GPAD}:
\begin{subequations}
	\begin{alignat}{5} 
	 &\bar \lambda_{ui} &&= \bar w_{ui} + \bar s_{ui},&&\quad\ \bar \lambda_{{\ell i}}= 0, \quad \forall i \in \II_{\bar u}\label{eq:yu=wu+su}\\
	 &\bar \lambda_{\ell i} &&= \bar w_{\ell i} + \bar s_{ \ell i},&&\quad\ \bar \lambda_{u i} = 0, \quad \forall i \in \II_{\bar \ell}\label{eq:yell=well+sell}\\
	 &\bar \lambda_{u,\JJ} &&\geq 0,&&\quad\ \bar \lambda_{\ell, \JJ} \geq 0 \label{eq:y>=0}
	\end{alignat}
    \label{eq:dual-conditions}
\end{subequations}
In other words, the \textrm{max} operator is not used for the elements $\bar \lambda_{{\ell},\II_{\bar \ell}}$, $\bar \lambda_{{u},\II_{\bar u}}$ (as they correspond to equality constraints) and the elements $\bar \lambda_{{\ell},\II_{\bar u}}$, $\bar \lambda_{{u},\II_{\bar \ell}}$ are just zeroed (as the corresponding constraints have been removed). 
Therefore, all QP relaxations have the same problem matrices and different executions of Algorithm~\ref{algo:GPAD} only differ in the projection step because of~\eqref{eq:dual-conditions}.
In particular, preconditioning can be computed just at the root node and maintained unaltered 
throughout the entire B\&B algorithm. We will further use~\eqref{eq:dual-conditions} for deriving the heuristic approach presented in Section~\ref{sec:heuristics_wobnb} for finding a suboptimal integer-feasible solution of the MIQP problem without using B\&B.

\subsection{Warm-starting the QP subproblems}

Conventionally the dual solution of the parent problem can be used as a feasible initial value for the children QP subproblems.
The following lemma provides a way of generating an initial guess, 
using not only the optimal dual solution of the parent problem but also the information of current branching variable $j$.
\begin{lemma}
	Let $z^*$, $\lambda^*$ be the primal-dual solution of the QP
	\begin{subequations}
		\begin{align}
		\displaystyle{\min_{z}}\quad
		& \frac{1}{2}z'Qz+c'z\\
		\st\quad            & \Aa_i z \leq b_i \quad \forall i \in {\II}\label{eq:QP1_wslemma}\\
	    & \Aa_i z = b_i \quad \forall i \in {\EE}\label{eq:QP2_wslemma}
		\end{align} \label{eq:QP_wslemma}
	\end{subequations}
	with $z \in \rr^{n}$, $\lambda \in \rr^{n_\lambda}$, and $\II \cup \EE = \{1,\ldots,n_\lambda\}$.
	Let $\II_{\textrm{na}}$ be the set of indices of inactive constraints at $z^*$,  
	\begin{align*}
	\Aa_{i} z^* < b_i, &\quad \forall i \in {\II_{\textrm{na}}} \\
	\Aa_{i} z^* = b_i, &\quad \forall i \in \{1,\ldots,n_\lambda\} \setminus {\II_{\textrm{na}}},                                                                        
	\end{align*}	                                                                                             where clearly $\II_{\textrm{na}} \subseteq \II$.
	                                                                                            
	For any $j \in \II_{\textrm{na}}$, consider the problem QP$_j$ obtained from~\eqref{eq:QP_wslemma} by moving $j$ from $\II$ to $\EE$, i.e., by using $\II \setminus \{j\}$ in~\eqref{eq:QP1_wslemma} and $\EE \cup \{j\}$ in~\eqref{eq:QP2_wslemma}. The dual vector $\bar \lambda$ defined by
    \begin{subequations}
	\begin{align}
	\bar{\lambda}_j &= \frac{\Aa_j z^* - b_j}{H_{jj}}\\
	\bar{\lambda}_i&=\lambda_i^*,~\forall i \in  \{1,\ldots,n_\lambda\}\setminus \{j\}
	\end{align}
    \label{eq:dual-warm-start}
    \end{subequations}
is feasible and such that the corresponding primal vector 
\begin{align}
\bar z = -Q^{-1}\Aa'\bar \lambda-Q^{-1}c     \label{eq:bar_z}
\end{align}
satisfies the equality constraint
\begin{align}
    \Aa_j\bar z - b_j=0
    \label{eq:Aj}
\end{align}
\label{lemma:warmstart}
\end{lemma}
\proof As $j\in \II_{\textrm{na}}$, we have that $\lambda^*_j=0$. Moreover $H_{jj}=\Aa_jQ^{-1}\Aa_j$.
Hence, from~\eqref{eq:bar_z} and the Karush-Kuhn-Tucker (KKT) conditions
of optimality of~\eqref{eq:QP_wslemma} $Qz^*+c+\Aa'\lambda^*=0$
we get
\beqarno
\Aa_j\bar z &=&-\sum_{i=1}^{n_\lambda}\Aa_jQ^{-1}\Aa_i\bar \lambda_i-\Aa_jQ^{-1}c\\
&=&-\sum_{i\neq j}\Aa_jQ^{-1}\Aa_i\lambda^*_i-\Aa_jQ^{-1}c-\Aa_jQ^{-1}\Aa_j\bar \lambda_j\\
&=&\Aa_jQ^{-1}\left(-c-\sum_{i=1}^{n_\lambda}\Aa_i\lambda^*_i\right)-(\Aa_jz^*-b_j)\\
&=&\Aa_jz^*-\Aa_jz^*+b_j=b_j
\eeqarno
\QED
Lemma~\ref{lemma:warmstart} allows us to obtain a good initial guess $\bar \lambda$ 
from the solution $z^*,\lambda^*$ of the parent node, in that the new equality
constraint~\eqref{eq:Aj} imposed in the child node is satisfied by construction.
Let the parent node be characterized by $\II_{\bar \ell}$, $\II_{\bar u}$
and optimal solution $z^*,\lambda^*$ and assume we are branching
on the $j$th integrality constraint, $j\in\{1,\ldots,p\}$. 
Since now on, we assume that the 
tuples $\TT_0=(\II_{\bar \ell}\cup{\{j\}},\II_{\bar u},\bar \lambda^0)$ 
and $\TT_1=(\II_{\bar \ell},\II_{\bar u}\cup{\{j\}},\bar \lambda^1)$ are pushed on the stack $\Ss$,
where
\begin{subequations}
 	\begin{align}
	\bar{\lambda}^0_{j_0} &= \frac{\bar{\ell}_j-\bar A_j z^*}{H_{j_0j_0}}\nonumber\\
	\bar{\lambda}^0_i&=\lambda_i^*,~\forall i \in  \{1,\ldots,n_\lambda\}\setminus \{j_0\}
	\end{align}
    and
	\begin{align}
	\bar{\lambda}_{j_1}^1 &= \frac{\bar A_j z^* - \bar u_j}{H_{j_1j_1}}\nonumber\\
	\bar{\lambda}_i^1&=\lambda_i^*,~\forall i \in  \{1,\ldots,n_\lambda\}\setminus \{j_1\}
	\end{align}
\label{eq:dual-warm-start2}
\end{subequations}
and $j_0$, $j_1$ are the indices corresponding to 
$\bar A_j z \geq \bar \ell_j$ and $\bar A_j z \leq \bar u_j$, respectively,
according to the ordering of the constraints defined in~\eqref{eq:df},
and $n_\lambda=2(m+p)$.

\subsection{Warm-starting binary variables} \label{sec:ws}

We now introduce a strategy in the MIQP solver that exploits warm starts on (all or some of) 
the binary constraints, that is an indication whether $\bar A_iz$ should be equal to $\bar\ell_i$ or $\bar u_i$.

Let $\tilde{\II}_{\bar \ell}$, $\tilde{\II}_{\bar u}$ be the sets containing the indices of the warm-started binary constraints equal to $\bar{\ell}_i$, $\bar{u}_i$ respectively, $\tilde{\II}_{\bar \ell}\cup\tilde{\II}_{\bar u}\neq\emptyset$,
$\tilde{\II}_{\bar \ell}\cap\tilde{\II}_{\bar u}=\emptyset$. Algorithm~\ref{algo:warmstart-binary} describes a new branching rule replacing Steps~\ref{step:branch-start}-\ref{step:branch-end} of Algorithm~\ref{algo:MIQP-GPAD}.
The combination of Algorithms~\ref{algo:MIQP-GPAD} and~\ref{algo:warmstart-binary} 
is executed as described below.

After solving the QP relaxation~\eqref{eq:QP_relaxation}
at the root node (steps up to Step~\ref{step:QP} of Algorithm~\ref{algo:MIQP-GPAD}), 
if $\bar A_iz^*_i\not\in\{\bar\ell_i,\bar u_i\}$ for some $i\in\{1,\ldots,p\}$, Algorithm~\ref{algo:warmstart-binary} is executed with $\JJ\leftarrow\{1,\ldots,p\}$
and proceeds to Step~\ref{step:branch-start_ws_bin}, where the variable with smallest index $j$ in either $\tilde{\II}_{\bar\ell}$ or $\tilde{\II}_{\bar u}$ is selected for branching. 
Two dual initial values are created at Step~\ref{step:lambda_ws_ws_bin} and two new subproblems are created at Step~\ref{step:AA_ws_bin}.

The Boolean variable \textsf{noQP} (Step~\ref{step:check-QP-to-be-solved}) keep tracks of whether the QP relaxation corresponding
to the new generated node must not be solved. This event happens when
the new index $j\in\tilde{\II}_{\bar\ell}$ or $j\in\tilde{\II}_{\bar u}$
and
when all indices in $\II_{\bar \ell}$ are contained in $\tilde{\II}_{\bar\ell}$,
all indices in $\II_{\bar u}$ are contained in $\tilde{\II}_{\bar u}$,
and indices in $\JJ$ still exist that belong to 
$\tilde{\II}_{\bar\ell}\cup\tilde{\II}_{\bar u}$. 
The latter condition is captured by measuring the cardinality $\card(\II_{\bar \ell})+1+\card(\II_{\bar u})$
of the new set obtained by adding $\{j\}$ to either $\II_{\bar \ell}$ or 
$\II_{\bar u}$ and verifying if it is smaller than or equal to the cardinality $\tilde{c}$ of 
$\tilde{\II}_{\bar\ell}\cup \tilde{\II}_{\bar u}$. 
In other words, if warm-start values still need
to be processed, the new node is marked as a node whose QP relaxation must not be solved.

Because of the last-in-first-out mechanism of $\Ss$, the subtree
originating from the node in which all warm-started indices are fixed to the corresponding
warm-start values gets processed first. If at some of the resulting leaf QP subproblem is feasible, the value of $V_0$ (the best known integer-feasible solution) is updated. The desired goal of the strategy of reducing the best known integer-feasible value $V_0$ as soon
as possible is then achieved, which in turns possibly reduces the overall number of QP relaxations required to find the global solution of the MIQP problem.

\begin{algorithm}[t]
	\caption{Warm-start for binary constraints}
	\label{algo:warmstart-binary}
	~~\textbf{Input}: Warm-start index sets $\tilde{\II}_{\bar \ell}$, $\tilde{\II}_{\bar u}$,
    $\tilde{c}\leftarrow\card(\tilde{\II}_{\bar \ell})+\card(\tilde{\II}_{\bar u})$;
    current sets $\JJ$, $\II_{\bar \ell}$, $\II_{\bar u}$, and optimal solution $z^*$, $\lambda^*$ of QP relaxation;
	\vspace*{.1cm}\hrule\vspace*{.1cm}
	\begin{enumerate}[label=\arabic*., ref=\theenumi{}]
		\item \label{step:branch-priorityto_ws_binary} \textbf{if} $\JJ \cap (\tilde{\II}_{\bar\ell} \cup \tilde{\II}_{\bar u})\neq \emptyset$ \textbf{then} 
		\begin{enumerate}[label=\theenumi{}.\arabic*., ref=\theenumi{}.\theenumii{}]
			\item \label{step:branch-start_ws_bin} $j \leftarrow \inf (\JJ \cap (\tilde{\II}_{\bar\ell} \cup \tilde{\II}_{\bar u}))$;  
			\item \label{step:lambda_ws_ws_bin} set $\bar \lambda^0$ and $\bar \lambda^1$ as in~\eqref{eq:dual-warm-start2}
            by activating $\bar A_j z=\bar \ell_j$ and $\bar A_j z=\bar u_j$, respectively;
            \item \label{step:AA_ws_bin}$\TT_0\leftarrow(\II_{\bar \ell}\cup\{j\},\II_{\bar u},\bar \lambda^0)$; $\TT_1\leftarrow(\II_{\bar \ell},\II_{\bar u}\cup\{j\},\bar \lambda^1)$; 
            \item \label{step:check-QP-to-be-solved}\textsf{noQP}$\leftarrow$ ($\II_{\bar \ell}\subseteq \tilde{\II}_{\bar\ell}$ \textbf{and}
            $\II_{\bar u}\subseteq \tilde{\II}_{\bar u}$ \textbf{and} $\card(\II_{\bar \ell})+1+\card(\II_{\bar u})\leq \tilde{c}$);			         
            \item \label{step:branch-end_ws_bin}\textbf{if} $j \in \tilde{\II}_{\bar \ell}$ \textbf{then}
    		\begin{enumerate}[label=\theenumii{}.\arabic*., ref=\theenumii{}.\theenumiii{}]
                \item push $\TT_1$ and then $\TT_0$ on $\Ss$;
                \item \label{step:QP-not-to-be-solved-1}\textbf{if} \textsf{noQP} \textbf{then} mark $\TT_0$ as a QP that must not be solved;
            \end{enumerate}
            \item \textbf{otherwise} 
            \begin{enumerate}[label=\theenumii{}.\arabic*., ref=\theenumii{}.\theenumiii{}]
                \item push $\TT_0$ and then $\TT_1$ on $\Ss$;
                \item \label{step:QP-not-to-be-solved-2}\textbf{if} \textsf{noQP} \textbf{then} mark $\TT_1$ as a QP that must not be solved; 
            \end{enumerate}
    	\end{enumerate}
		\item \textbf{otherwise} 
        \begin{enumerate}[label=\theenumi{}.\arabic*., ref=\theenumi{}.\theenumii{}]
            \item \label{step:regular-branching} execute Step~\ref{step:check-integrality} of Algorithm~\ref{algo:MIQP-GPAD};
    	\end{enumerate}
	\end{enumerate}
	\vspace*{.1cm}\hrule\vspace*{.1cm}
	~~\textbf{Output}: Two new tuples $\TT_1,\TT_0$ (in appropriate order) are pushed on stack $\Ss$.
\end{algorithm}

Note that at Step~\ref{step:regular-branching} of Algorithm~\ref{algo:warmstart-binary}
to warm start also real variables when running Step~\ref{step:check-integrality} of Algorithm~\ref{algo:MIQP-GPAD} we need to set $\TT_0\leftarrow(\II_{\bar \ell}\cup\{j\},\II_{\bar u},\bar \lambda^0)$; $\TT_1\leftarrow(\II_{\bar \ell},\II_{\bar u}\cup\{j\},\bar \lambda^1)$ and $\bar \lambda^0$ and $\bar \lambda^1$ defined as in~\eqref{eq:dual-warm-start2}. Note also that at Steps~\ref{step:QP-not-to-be-solved-1},~\ref{step:QP-not-to-be-solved-2} of Algorithm~\ref{algo:warmstart-binary}
the QP problem is marked not to be solved. When the corresponding node
will be popped up from $\Ss$, the warm-start value $\bar \lambda$ of real variables
will not be updated, but will remain equal to the dual solution $\lambda^*$ of the last ancestor (QP relaxation) which was solved.

\begin{figure}[h!]
	\centering
	\begin{forest}
		for tree={ align=center,draw,fill=black!5,circle,minimum size=1.8em,inner sep=0pt,text height=1.7ex, text depth=0.2ex}
		[1,label={$\star$$\star$$\star$}
		[7,label=$0$$\star$$\star$
		[,label=$00\star$
		[,label=below:$000$]
		[,label=below:$001$]
		]
		[,label=$01\star$
		[,label=below:$010$]
		[,label=below:$011$]
		]
		]
		[,label=$1$$\star$$\star$,fill=red!25,thick,dashed
		[,label=$10\star$,fill=red!25,thick,dashed
		[2,label=below:$100$,fill=red!25]
		[3,label=below:$101$,fill=red!25]
		]
		[4,label=$11\star$
		[5,label=below:$110$]
		[6,label=below:$111$]
		]
		]
		]
	\end{forest}
	\caption{Illustration example with 3 binary variables and $(1, 0, \star)$ as a binary warm start. The numbers denote the order in which the QP relaxations are solved, dashed nodes correspond to QP relaxations that are not solved.}
	\label{fig:tree_exploration_ws_bin}
\end{figure}
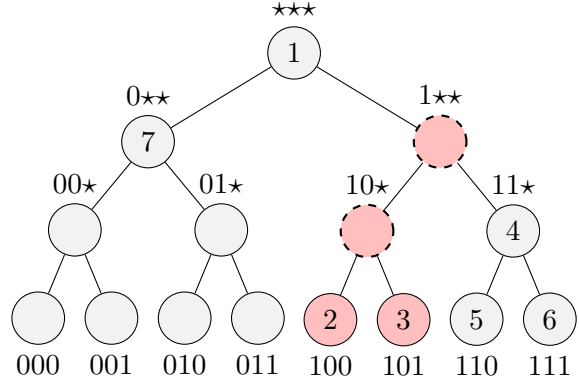

We illustrate the tree generated by Algorithms~\ref{algo:MIQP-GPAD},~\ref{algo:warmstart-binary} with
an example with 3 binary constraints $z_1,z_2,z_3\in\{0,1\}$, warm started using the sequence $(1, 0, \star)$,
i.e., only $z_1$, $z_2$ are warm started at $z_1=1$, $z_2=0$. The tree is shown in Fig.~\ref{fig:tree_exploration_ws_bin},
where the number inside the node denotes the order in which the QP relaxations are executed.

First the root node is solved and then two problems are pushed on the stack, with node $(1, \star, \star)$ 
on top of stack $\Ss$. Due to the provided warm start, the leaf problem $(1, 0, 0)$ gets solved first (QP \#2), and 
$(1, 0, 1)$ immediately after (QP \#3). Next, problem $(1, \star, \star)$ is popped from $\Ss$ (due to the last-in-first-out method), but not solved (depicted as a dashed circle in Fig.~\ref{fig:tree_exploration_ws_bin}), as it was marked like that
at Step~\ref{step:QP-not-to-be-solved-2} of Algorithm~\ref{algo:warmstart-binary}. Similarly, problem $(1, 0, \star)$ is popped but not solved. 
Next, problem  $(1, 1, \star)$ is popped from $\Ss$ and solved (QP \#4).

Note that we have saved solving 2 QP relaxations without compromising the optimality of the MIQP solution. In general, using this approach at least $\card(\tilde{\II}_{\bar \ell})+\card(\tilde{\II}_{\bar u})$ QP relaxations can be saved.

\begin{figure}[h!]
	\centering
	\begin{forest}
		for tree={ align=center,draw,fill=black!5,circle,minimum size=1.8em,inner sep=0pt,text height=1.7ex, text depth=0.2ex}
		[1,label={$\star$$\star$$\star$}
		[,label=$\star$$\star$$0$,fill=red!25,thick,dashed
		[,label=$\star00$,fill=red!25,thick,dashed
		[2,label=below:$000$,fill=red!25]
		[4,label=below:$100$]
		]
		[5,label=$\star10$
		[6,label=below:$010$]
		[7,label=below:$110$]
		]
		]
		[,label=$\star$$\star$$1$,fill=red!25,thick,dashed
		[,label=$\star01$,fill=red!25,thick,dashed
		[3,label=below:$001$,fill=red!25]
		[,label=below:$101$]
		]
		[,label=$\star11$
		[,label=below:$010$]
		[,label=below:$111$]
		]
		]
		]
	\end{forest}
	\caption{Illustration example with 3 binary variables and $(0, 0, \star)$ as a binary warm start. The numbers denote the order in which the QP relaxations are solved, dashed nodes correspond to QP subproblems that are ignored.}
	\label{fig:tree_exploration_ws_bin_2}
\end{figure}
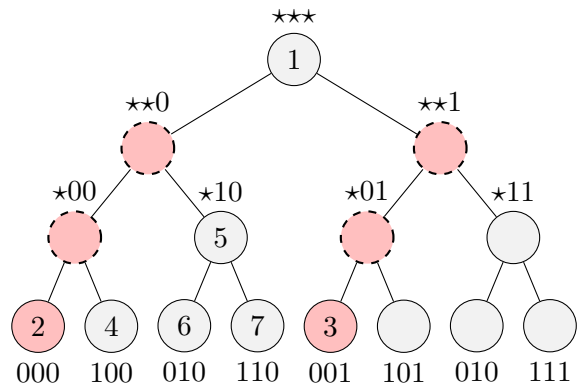

A few more QPs may be saved if different branching rules are used. Consider the example shown in Fig.~\ref{fig:tree_exploration_ws_bin_2}, in which the binary constraint to branch 
is selected according to the maximum fractional part in Step~\ref{step:branch-start_ws_bin} of Algorithm~\ref{algo:warmstart-binary}, the binary warm start is $(0, 0, \star)$,
and Step~\ref{step:check-QP-to-be-solved} is changed to
\begin{itemize}
\item []
\textsf{noQP}$\leftarrow$ ($\II_{\bar \ell}\cap\tilde{\II}_{\bar u}=\emptyset$ \textbf{and}
$\II_{\bar u}\cap \tilde{\II}_{\bar \ell}=\emptyset$ \textbf{and} $\card(\II_{\bar \ell}\cap \tilde{\II}_{\bar \ell})+1+\card(\II_{\bar u}\cap\tilde{\II}_{\bar u})<\tilde{c}$)
\end{itemize}

In this case, as shown in Fig.~\ref{fig:tree_exploration_ws_bin_2},
we save 4 QP subproblems. We remark that saving of more QP relaxation may vary depending upon the problem under consideration, specifically on how the tree layout progresses (due to branching) with respect to the provided warm-start.

\subsubsection{Warm-starting binary variables correspond to an exclusive-or condition}\label{sec:ws_sos1}

We consider a specific subclass of MIQP problems where the binary variables are subject to an exclusive-or condition, such as due to one-hot encoding of categorical decision variables. Such problems lies in category of the Special Ordered Sets of Type 1 (SOS1 or S1), defined as a set of variables for which only one binary variable can take a nonzero value while the remaining are zero.
For such problems, we extend the idea of warm-start with~\textsf{noQP} for efficiently reducing the number of sub-problems to be explored. This is done by replacing Step~\ref{step:regular-branching} of Algorithm~\ref{algo:warmstart-binary} with Algorithm~\ref{algo:warmstart-sos1}, where $\tilde{s}$ denotes the cardinality of the binaries correspond to the submodels, also referred to as modes. We remark that for simplicity of notation, we describe the idea in Algorithm~\ref{algo:warmstart-sos1} for the case where all the binaries are subjected to an exclusive-or condition i.e. $\tilde{s}=p$. 

\begin{algorithm}[htb!]
	\caption{\textsf{noQP} for binaries corresponding to submodels}\label{algo:warmstart-sos1}
	\begin{enumerate}[label=\arabic*., ref=\theenumi{}]
		\item \textbf{if} $\JJ=\emptyset$ 
		\textbf{or} 
		$t_i\eqdef\bar A_iz^*\in\{\bar \ell_i,\bar u_i\}$, 
		$\forall i\in \JJ$ 
		\textbf{then}
		$V_0\leftarrow V^*$, $\zeta^*\leftarrow z^*$; \textbf{otherwise}
	\begin{enumerate}[label=\theenumi{}.\arabic*., ref=\theenumi{}.\theenumii{}]
		\item $\displaystyle{j\leftarrow\arg\min_{i\in
				J}}$ $\left|t_i-\frac{\bar \ell_i+ \bar u_i}{2}\right|$;
		\item set $\bar \lambda^0$ and $\bar \lambda^1$ as in~\eqref{eq:dual-warm-start2}
		by activating $\bar A_j z=\bar \ell_j$ and $\bar A_j z=\bar u_j$, respectively;
		\item $\TT_0\leftarrow(\II_{\bar \ell}\cup\{j\},\II_{\bar u},\bar \lambda^0)$; $\TT_1\leftarrow(\II_{\bar \ell},\II_{\bar u}\cup\{j\},\bar \lambda^1)$; 
	\item \textbf{if} $t_{j}\leq\frac{\bar \ell_i+\bar u_i}{2}$ \textbf{then}
	\begin{enumerate}[label=\theenumii{}.\arabic*., ref=\theenumii{}.\theenumiii{}]
		\item push $\TT_1$ and then $\TT_0$ on $\Ss$;
		\item \textsf{noQP}$_{\TT_0}$$\leftarrow$ ($\card(\II_{\bar \ell})+1= \tilde{s}$ \textbf{or} $\card(\II_{\bar u})> 1$ \textbf{or} $\card(\II_{\bar \ell})+2+\card(\II_{\bar u})= \tilde{s}$);
		\item \textbf{if} \textsf{noQP}$_{\TT_0}$ \textbf{then} mark $\TT_0$ as a QP that must not be solved;
	\end{enumerate}
	\item \textbf{otherwise} 
	\begin{enumerate}[label=\theenumii{}.\arabic*., ref=\theenumii{}.\theenumiii{}]
		\item push $\TT_0$ and then $\TT_1$ on $\Ss$;
		\item \textsf{noQP}$_{\TT_1}$$\leftarrow$ ($\card(\II_{\bar u})\geq 1$ \textbf{or} $\card(\II_{\bar \ell})+2+\card(\II_{\bar u})= \tilde{s}$);
		\item \textbf{if} \textsf{noQP}$_{\TT_1}$ \textbf{then} mark $\TT_1$ as a QP that must not be solved; 
	\end{enumerate}
\end{enumerate}
	\end{enumerate}
\end{algorithm}

\begin{figure*}[thb!]
	\centering
	\begin{forest}
		/tikz/every node/.append style={font=\small},
		for tree={ align=center,draw,fill=black!5,circle,minimum size=1.8em,inner sep=0pt,text height=1.7ex, text depth=0.2ex}
		[1,label={$\star$$\star$$\star$$\star$}
		[,label=$0$$\star$$\star$$\star$,fill=red!45,thick,dashed
		[,label=$00$$\star$$\star$,fill=red!55,thick,dashed
		[,label=$000\star$,fill=red!10,thick,dashed
		[,label=below:$0000$,fill=red!10,thick,dashed]
		[,label=below:$0001$,fill=red!10,thick,dashed]]
		[,label=$001\star$,fill=red!10,thick,dashed
		[,label=below:$0010$,fill=red!10,thick,dashed]
		[,label=below:$0011$,fill=red!10,thick,dashed]
		]
		]
		[,label=$01$$\star$$\star$,fill=red!45,thick,dashed
		[,label=$010\star$,fill=red!45,thick,dashed
		[,label=below:$0100$,fill=red!10,thick,dashed]
		[2,label=below:$0101$,fill=red!45]
		]
		[,label=$011\star$,fill=red!45,thick,dashed
		[3,label=below:$0110$,fill=red!45]
		[,label=below:$0111$,fill=red!10,thick,dashed]
		]
		]
		]
		[4,label=$1$$\star$$\star$$\star$
		[,label=$10$$\star$$\star$
		[,label=$100\star$,
		[,label=below:$1000$,fill=red!10,thick,dashed]
		[,label=below:$1001$]
		]
		[,label=$101\star$
		[,label=below:$1010$]
		[,label=below:$1011$,fill=red!10,thick,dashed]
		]
		]
		[,label=$11$$\star$$\star$,fill=red!10,thick,dashed
		[,label=$110\star$,fill=red!10,thick,dashed
		[,label=below:$1100$,fill=red!10,thick,dashed]
		[,label=below:$1101$,fill=red!10,thick,dashed]
		]
		[,label=$111\star$,fill=red!10,thick,dashed
		[,label=below:$1110$,fill=red!10,thick,dashed]
		[,label=below:$1111$,fill=red!10,thick,dashed]
		]
		]
		]
		]
	\end{forest}
	\caption{Illustrative example with 4 binary variables with constrains $z_1+z_2=1$ and $z_3+z_4=1$ and $(0, 1, \star, \star)$ as a binary warm start. The numbers denote the order in which the QP relaxations are solved, 
	dark-dashed nodes correspond to QP relaxations that are explored but not solved.
		Light-dashed nodes correspond to the infeasible QP subproblems not to be solved.}
	\label{fig:sos1_ex}
\end{figure*}
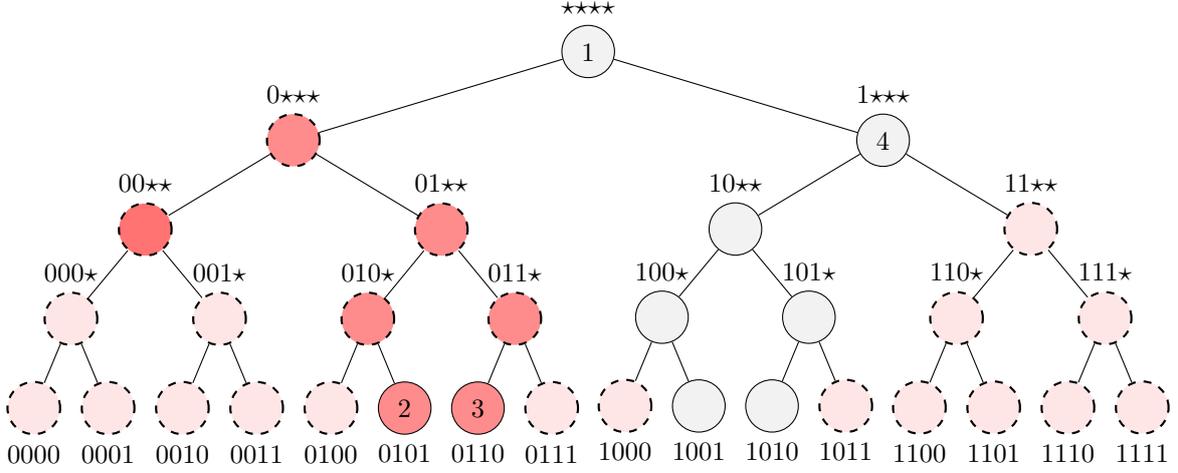

We illustrate the tree generated by Algorithms~\ref{algo:MIQP-GPAD},~\ref{algo:warmstart-binary},~\ref{algo:warmstart-sos1} with
an example with 4 binary variables having $\tilde{s}=2$, constraints  $z_1+z_2=1$ and $z_3+z_4=1$, warm started using the sequence $(0, 1, \star, \star)$,
i.e., $z_1$, $z_2$ are warm started with $z_1=0$, $z_2=1$. The tree is shown in Fig.~\ref{fig:sos1_ex}, the numbers denote the order in which the QP relaxations are solved, dark-dashed nodes correspond to QP relaxations that are explored but not solved.
Light-dashed nodes correspond to the infeasible QP subproblems not to be solved.

\section{Heuristic methods for suboptimal binary-feasible MIQP solutions}\label{sec:heuristics} 

In this section, we propose two
heuristic approaches to possibly obtain suboptimal integer feasible solutions of the  MIQP problem~\eqref{eq:QP}. The first approach does not use B\&B, the second one only branches on a reduced
set of binary constraints. Both approaches aim at finding an upper bound on the cost of the MIQP problem 
upfront before executing the full B\&B algorithm.
Indeed, denoting by $V_{H}^{*}$ the suboptimal integer feasible solution of the heuristic approach,
$V_{H}^{*}\in\rr\cup\{+\infty\}$, QP relaxations in which the cost goes above $V_{H}^{*}$ can be terminated prematurely by checking condition~\eqref{eq:V0} and the corresponding node is fathomed. Moreover, 
if $V_{H}^{*}<+\infty$ and B\&B fails to find a better integer feasible solution within a specified time limit, 
or B\&B is not executed at all, a suboptimal solution is available.

\subsection{Heuristic approach without using B\&B}\label{sec:heuristics_wobnb}

In the sequel we denote by $z^{*}_{H}$ the solution of the presented heuristic approach, if any
is found. Since either $\bar A_iz^{*}_{H}=\bar\ell_i$ or $\bar A_iz^{*}_{H}=\bar u_i$ must hold to satisfy the integrality constraints~\eqref{eq:QP4}, for all $i=1,\ldots,p$ the corresponding dual variables are such that either $\bar \lambda_{ u_{i}H}^{*}=0$
and $\bar \lambda_{\ell_{i}H}^{*}$ unconstrained/unrestricted in sign, or vice versa (similar to~\eqref{eq:yu=wu+su}-\eqref{eq:yell=well+sell}). In other words, 
the vector $\smallmat{\bar \lambda_{ u_{i}H}^{*}\\\bar \lambda_{\ell_{i}H}^{*}}$ must belong to 
the nonconvex set given by the union of the orthogonal real axes.
We propose the heuristic method described in Algorithm~\ref{algo:heuristic} to
define the values of $\bar \lambda_{\ell}$, $\bar \lambda_{u}$ during the execution
of Algorithm~\ref{algo:GPAD} to impose such a nonconvex constraint,
analogously to the approach described in~\cite{TMBB17} for ADMM.

\begin{algorithm}[b]
	\caption{Heuristic for suboptimal binary feasible MIQP solution}\label{algo:heuristic}
	\begin{enumerate}[label=\arabic*., ref=\theenumi{}]
        \item $\smallmat{{\lambda_{u_{k+1}}}\\{\lambda_{\ell_{k+1}}}}\leftarrow \mathrm{max}\left\lbrace  \smallmat{{w_{u_k}}+{s_{u_k}}\\{w_{\ell_k}}+{s_{\ell_k}}},0\right\rbrace $;
		\item \label{step:fori-in-1p}\textbf{for} $i=1,\ldots,p$ \textbf{do}
        \begin{enumerate}[label=\theenumi{}.\arabic*., ref=\theenumi{}.\theenumii{}]
            \item \textbf{if} $\bar{A}_{i}z_{k} \geq \tfrac{\bar l_{i}+\bar u_{i}}{2}$ \label{step:heu condition Aiz>=l+u/2} 
            \textbf{then} $\bar \lambda_{u_{{i}_{k+1}}} \leftarrow\bar w_{ u_{{i}_{k}}}+\bar s_{u_{{i}_{k}}}$;
	       	$\bar \lambda_{{\ell}_{{i}_{k+1}}}\leftarrow0$;
        \end{enumerate}
        \item [] \textbf{otherwise}
        \begin{enumerate}[label=\theenumi{}.\arabic*., ref=\theenumi{}.\theenumii{}]
            \setcounter{enumii}{1}
            \item 
    		$\bar \lambda_{{u}_{{i}_{k+1}}}\leftarrow0$;
    		$\bar \lambda_{ \ell_{{i}_{k+1}}} \leftarrow \bar w_{ \ell_{{i}_{k}}}+\bar s_{ \ell_{{i}_{k}}}$;
        \end{enumerate}
    \end{enumerate}
\end{algorithm}

Assuming that Algorithm~\ref{algo:heuristic} converges
under projection on such a nonconvex constraint set, 
$z^{*}_{H}$ satisfies~\eqref{eq:QP4} by construction,
as all dual variables corresponding to the constraints~\eqref{eq:QP4}
are either not restricted in sign, which is equivalent to treating
the corresponding constraint as an equality constraint, or set equal to zero,
which is equivalent to ignoring the corresponding constraint (similar to steps~\eqref{eq:yu=wu+su}-\eqref{eq:yell=well+sell}). We propose to
apply Algorithm~\ref{algo:heuristic} only
after the QP relaxation~\eqref{eq:QP_relaxation}
with $\II_{\bar u}=\II_{\bar \ell}=\emptyset$
is solved by Algorithm~\ref{algo:GPAD} to optimality and the problem
is found feasible but not integer feasible. In this case, 
Algorithm~\ref{algo:GPAD} is executed again from Step~\ref{step:repeat}, where now Step~\ref{step:yk1} is replaced by Algorithm~\ref{algo:heuristic}, until stopping criteria~\eqref{eq:gpad feas cond} and~\eqref{eq:gpad opt cond} are satisfied.
We replace $-\smallmat{w_{k}\\w_{eq,k}}'\smallmat{s_{k}\\s_{eq,k}}$ with $|-\smallmat{w_{k}\\w_{eq,k}}'\smallmat{s_{k}\\s_{eq,k}}|$ and exclude $w_{k}\geq0$
in condition~\eqref{eq:gpad opt cond} when executing Algorithm~\ref{algo:GPAD}.

\subsection{Mid-way heuristic approach}\label{sec:midway_heuristic}

We propose a second approach where we fix a subset of binary constraints to either $\bar{A}_{i}z=\bar\ell_i$ or $\bar{A}_{i}z=\bar u_i$ using a heuristic and then treat the remaining ones using standard B\&B.
For this reason, we denote the method as \emph{mid-way heuristic approach}.

This idea is summarized by Algorithm~\ref{algo: a midway approach}, where the sets  $\mathcal K_{\bar \ell}, \mathcal K_{\bar u}$ and $\LL$ are defined which corresponds to the binary constraints satisfying $\bar A_{i}z^{*}_{R} \leq \epsilon_{\bar \ell}$,  $\bar A_{i}z^{*}_{R} \geq \epsilon_{\bar u}$ and $\epsilon_{\bar \ell} < \bar A_{i}z^{*}_{R} < \epsilon_{\bar u}$, respectively, with $\KK_{\bar u} \cap \KK_{\bar\ell}=
\KK_{\bar u} \cap \LL = \KK_{\bar\ell} \cap \LL=\emptyset$, $\KK_{\bar u} \cup \KK_{\bar \ell} \cup \LL= \{1, \ldots, p\}$ where $\bar\ell_i \leq \epsilon_{\bar \ell} < \epsilon_{\bar u} \leq \bar u_i$. 

\begin{algorithm}[!thb]
	\caption{the mid-way heuristic approach}\label{algo: a midway approach}
	~~\textbf{Input}: solution $z^{*}_{R}$ 
    of QP relaxation~\eqref{eq:QP_relaxation} with $\II_{\bar u}=\II_{\bar \ell}=\emptyset$, thresholds $\epsilon_{\bar \ell}$, $\epsilon_{\bar u}$.
	\vspace*{.1cm}\hrule\vspace*{.1cm}
		\begin{enumerate}[label=\arabic*., ref=\theenumi{}]
		\item \label{step: vio prone init}$\mathcal K_{\bar \ell}, \mathcal K_{\bar u}, \LL\leftarrow \emptyset$
		\item $\mathcal K_{\bar \ell}\leftarrow \{i_{\ell}\in \{1,\ldots,p\} | \bar{A}_{i_{\ell}}z^{*}_{R} \leq \epsilon_{\bar \ell}\}$
		\item $\mathcal K_{\bar u}\leftarrow \{i_{u}\in \{1,\ldots,p\} | \bar{A}_{i_{u}}z^{*}_{R} \geq \epsilon_{\bar u}\}$
		\item $\mathcal L \leftarrow \{1, \ldots, p\} \setminus \{\mathcal K_{\bar \ell}\cup \mathcal K_{\bar u} \}$	
	\end{enumerate}
\vspace*{.1cm}\hrule\vspace*{.1cm}
~~\textbf{Output}: sets  $\mathcal K_{\bar \ell}, \mathcal K_{\bar u}, \LL$.
\end{algorithm}

First, the QP relaxation~\eqref{eq:QP_relaxation}
with $\II_{\bar u}=\II_{\bar \ell}=\emptyset$
is solved by Algorithm~\ref{algo:GPAD} to optimality and if the problem
is found feasible but not integer feasible, then Algorithm~\ref{algo: a midway approach} is executed. Consequently after the sets $\mathcal K_{\bar \ell}$, $\mathcal K_{\bar u}$, $\LL$ are generated,
the heuristic approach of Algorithm~\ref{algo:heuristic} is run as described in Section~\ref{sec:heuristics_wobnb}.
This is followed by the B\&B Algorithm~\ref{algo:MIQP-GPAD}, where
the binary variables belonging to the sets $\KK_{\bar \ell}, \KK_{\bar u}$ are fixed to $\bar{A}_{i}z=\bar{\ell}_{i}$, $\bar{A}_{i}z = \bar{u}_{i}$ respectively.
If a suboptimal integer feasible solution $V^*_H$ is returned
by Algorithm~\ref{algo:heuristic}, an upper bound on the cost as in Step~\ref{step:init} 
is available for B\&B Algorithm~\ref{algo:MIQP-GPAD}. The latter is initialized with 
$V_0\leftarrow V^*_H$, run directly from Step~\ref{step:branch-start} (the QP relaxation
at the root node has already been solved to compute $z^*_R$), 
and allowed to only branch on indices $i\in\LL$.

Note that in applications where throughput is limited, one can control the amounts of computations
required by B\&B by calibrating the thresholds $\epsilon_{\bar \ell}$, $\epsilon_{\bar u}$.
The effectiveness of the heuristic will be analyzed in Section~\ref{sec:results}.

\subsection{Mid-way approach for binary warm start}

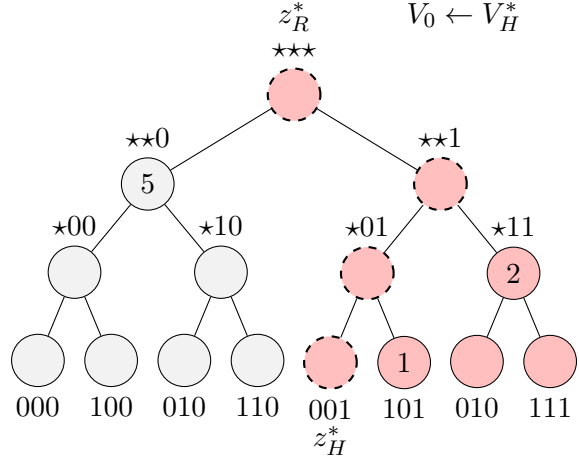
\begin{figure}[h!]
	\centering
	\begin{forest}
		for tree={ align=center,draw,fill=black!5,circle,minimum size=1.8em,inner sep=0pt,text height=1.7ex, text depth=0.2ex}
		[,label={$\star$$\star$$\star$},fill=red!25,thick,dashed,name = zR  
		[5,label=$\star$$\star$$0$
		[,label=$\star00$
		[,label=below:$000$]
		[,label=below:$100$]
		]
		[,label=$\star10$
		[,label=below:$010$]
		[,label=below:$110$]
		]
		]
		[,label=$\star$$\star$$1$,fill=red!25,thick,dashed
		[,label=$\star01$,fill=red!25,thick,dashed
		[,label=below:$001$,fill=red!25,thick,dashed, name=zH] 
		[1,label=below:$101$,fill=red!25]
		]
		[2,label=$\star11$,fill=red!25
		[,label=below:$010$,fill=red!25]
		[,label=below:$111$,fill=red!25]
		]
		]
		]
		\node [below=of zH, anchor=south] {$z^*_H$};
		\node [above=of zR, anchor=north, name=A] {$z^*_R$};
		\node [right =of A] {$V_0\leftarrow V^*_H$};
	\end{forest}
	\caption{Illustration example with 3 binary variables and $z^*_H=(0,0,1)$ as solution of heuristic approach, $\mathcal K_{\bar \ell}=\emptyset,~\mathcal K_{\bar u}=\{3\},~\LL=\{1,2\}$ resulting in $(\star, \star, 1)$ as a binary warm start. The numbers denote the order in which the QP relaxations are solved, 
	dashed nodes correspond to QP subproblems that are ignored. The red nodes shows the prioritized sub-tree.
	}
	\label{fig:tree_exploration_ws_bin_2_midway}
\end{figure}

The sets $\KK_{\bar \ell}$, $\KK_{\bar u}$, and $\LL$ generated
by Algorithm~\ref{algo: a midway approach} can be used for binary warm start in Algorithm~\ref{algo:warmstart-binary} 
by setting
\begin{align}
&\tilde{\II}_{\bar \ell} \leftarrow\KK_{\bar \ell},\ \tilde{\II}_{\bar u}\leftarrow\KK_{\bar u} 
\end{align}
Consider an example of MIQP with 3 binary variables, let $z^*_H=(0,0,1)$ be the solution provided
by Algorithm~\ref{algo:heuristic}, and let $\mathcal K_{\bar \ell}=\emptyset$, $\mathcal K_{\bar u}=\{3\}$,
$\LL=\{1,2\}$. 
 This leads to the warm start $z_3=1$ (as $\tilde{\II}_{\bar u}\leftarrow\KK_{\bar u}$) 
as shown in Fig.~\ref{fig:tree_exploration_ws_bin_2_midway}.

First, the best known integer feasible solution $V_0$ is initialized with $V^*_H$. From the extracted solution $z^*_R$ of the QP relaxation at the root node ($\star$, $\star$, $\star$), two problems are pushed on the stack with high priority to $(\star,\star,1)$ (due to the provided warm-start $z_3=1$), which is marked
as not to be solved.
 
After popping the tuple associated with $(\star,\star,1)$,
two problems are pushed on the stack with high priority to the $(\star,0,1)$ node, which then is popped from the stack (due to the last-in-first-out method), but again not solved (depicted as a dashed circle in Fig.~\ref{fig:tree_exploration_ws_bin_2_midway}). The leaf problem $(1,0,1)$ is popped from the stack and solved (QP \#1). The leaf problem $(0,0,1)$ is not solved again, as its solution $z^*_H$ is already
available. Then problem $(\star,1,1)$ is popped from the stack and solved (QP \#2), and standard B\&B 
proceeds on the remaining nodes.

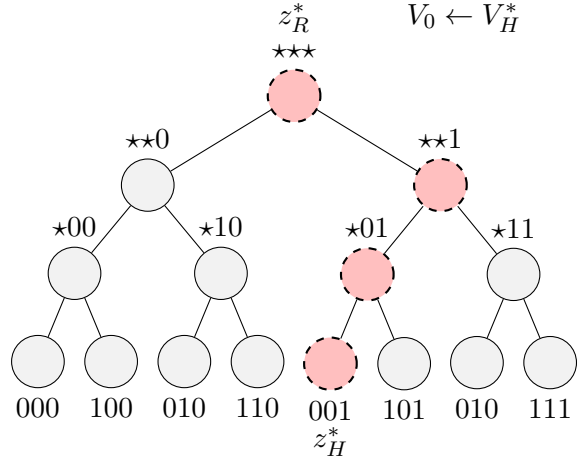
\begin{figure}[tb!]
	\centering
	\begin{forest}
		for tree={ align=center,draw,fill=black!5,circle,minimum size=1.8em,inner sep=0pt,text height=1.7ex, text depth=0.2ex}
		[,label={$\star$$\star$$\star$},fill=red!25,thick,dashed,name = zR  
		[,label=$\star$$\star$$0$
		[,label=$\star00$
		[,label=below:$000$]
		[,label=below:$100$]
		]
		[,label=$\star10$
		[,label=below:$010$]
		[,label=below:$110$]
		]
		]
		[,label=$\star$$\star$$1$,fill=red!25,thick,dashed
		[,label=$\star01$,fill=red!25,thick,dashed
		[,label=below:$001$,fill=red!25,thick,dashed, name=zH] 
		[,label=below:$101$]
		]
		[,label=$\star11$
		[,label=below:$010$]
		[,label=below:$111$]
		]
		]
		]
		\node [below=of zH, anchor=south] {$z^*_H$};
		\node [above=of zR, anchor=north, name=A] {$z^*_R$};
		\node [right =of A] {$V_0\leftarrow V^*_H$};
	\end{forest}
	\caption{Illustration example with 3 binary variables and $z^*_H=(0,0,1)$ as solution of heuristic approach. The dashed nodes correspond to QP subproblems that are ignored.}
	\label{fig:tree_exploration_ws_bin_2_heu}
\end{figure}

Similarly, the heuristic approach presented in Section~\ref{sec:heuristics_wobnb} can be combined with B\&B and binary warm start framework. 
Fig.~\ref{fig:tree_exploration_ws_bin_2_heu} illustrates this idea using an example of MIQP with 3 binary variables having $z^*_H=(0,0,1)$ as the solution found using Algorithm~\ref{algo:heuristic}. 
Hence, in general this proposed framework can be used when any kind of presolving techniques are combined with B\&B. 

\section{Numerical results}
\label{sec:results}

The numerical experiments presented in this section were performed on a desktop computer with Intel Core i7-4700MQ CPU with 2.40 GHz and 8 GB of RAM, running MATLAB R2015a.
Algorithms~\ref{algo:MIQP-GPAD},~\ref{algo:warmstart-binary},~\ref{algo:warmstart-sos1} were implemented in interpreted MATLAB code and Algorithms~\ref{algo:GPAD},~\ref{algo:Infeasibility detection},~\ref{algo:heuristic},~\ref{algo: a midway approach} in Embedded MATLAB code and compiled. Identical tolerance values are used in GUROBI for a fair comparison.

\subsection{Branch \& Bound on random MIQPs}\label{sec:exp_bnb_random}

The B\&B method of Algorithm~\ref{algo:MIQP-GPAD}, denoted as \textsf{miqpGPAD}, is tested on randomly generated MIQP problems with condition number $\kappa=10$ of the primal Hessian $Q$. We set
$Q=U\Sigma V'$, where $U,V$ are orthogonal matrices generated by QR decomposition of random $n\times n$ matrices, and $\Sigma$ is diagonal with nonzero entries having logarithms equally spaced between $\pm\log(\kappa)/4$~\cite{BTT91}. The entries of matrix $A$ are generated from the normal distribution $\NN(0,0.0025)$,
	$\ell$, $u$ from the uniform distribution $\UU(0,100)$, $c$ from $\NN(0,1)$. 

\begin{table}[h!]
	\centerline{\small\begin{tabular}{*6r}
			\hline
			$n$ & $m$ & $p$&$q$ &\textsf{miqpGPAD} & GUROBI \\
			\hline
			10 &  100 &   2 &2 &  6.11 & 3.52  \\
			50 &   25 &   5 &3 &  2.59 & 8.07  \\
			50 &  150 &  10 &5 & 23.16 & 31.47 \\
			100 &   50 &   2&5 &  4.72 & 25.62 \\
			100 &  200 &  15&5 & 47.48 & 111.04\\
			150 &  100 &   5&5 & 10.64 & 70.35 \\
			150 &  200 &  20&5 & 70.18 & 192.23 \\
			200&   50  & 15	&6 & 12.92 & 116.41 \\
			\hline
		\end{tabular}
	}\vspace*{.5cm}
	\caption{Average CPU time (ms) on random MIQP problems over 50 instances for each combination of $n$, $m$, $p$, $q$.}
	\label{tab:random-MIQP-comparison}
\end{table}
The tolerance values considered in Algorithm~\ref{algo:GPAD} are $ \epsilon_{V},\epsilon_{G}=1e-5$, while in Algorithm~\ref{algo:Infeasibility detection} we set $\epsilon_{I}=1e-2$. The CPU time reported for solving feasible MIQP problems, averaged over 50 executions, is listed in Table~\ref{tab:random-MIQP-comparison}, where we report the number $n$ variables, $m$ of inequality constraints, $p$ of binary constraints, and $q$ of equality constraints. 
The results show that the computation time of proposed scheme is comparable to that of the commercial solver GUROBI on such relatively small-size and well-conditioned MIQPs.

\subsection{Piecewise-affine regression problem}

Next, we consider a PieceWise Affine (PWA) regression problem that we solve using the MIQP formulation of~\cite{NMPB17MED}. This requires solving a MIQP at every sampling instant in a receding-horizon fashion. 
We consider data generated by the single-input and single-output (SISO)  {PWA AutoRegressive with eXogenous input} (PWARX) system from~\cite{MNPB_pwarx_ijrnc}
\begin{eqnarray*}
y(k)&=& \left\{
    \begin{array}{ll}
A_1\smallmat{1\\x(k)}+ e_{\mathrm{o}}(k) &  \text{if } -0.3x_{1}(k)+0.6x_{3}(k) +0.3 > 0\\
A_2\smallmat{1\\x(k)}+ e_{\mathrm{o}}(k) &  \text{if } -0.3x_{1}(k)+0.6x_{3}(k) +0.3 < 0
    \end{array}
\right.5\\
&&A_1={\begin{bmatrix}0.2 & 0.5 &-0.1 &1 &0.2 \end{bmatrix}}\\
&&A_2={\begin{bmatrix}  -0.3 & 0.8  & 0.1  & 0.4 & 0.05\end{bmatrix}}
\end{eqnarray*}
where $e_{\mathrm{o}}$ is a zero-mean white Gaussian process with variance $0.01$, leading to a {Signal-to-Noise Ratio} (SNR) of $22$ dB. 

The  regressor vector is $x(k)=[y(k-1)\ y(k-2)\ u(k-1)\ u(k-2)]'$,
where $u(k)$ is a randomly  generated external input  from a uniform  distribution having values between $-2$ and $2$. 
The training dataset consists of $N=1000$ samples, the validation dataset of $N_\mathrm{val}=300$ samples. 
A PWA model  with $\tilde{s}=2$ modes is considered  with horizon length $T=2$ steps. 
Each formulated MIQP sub-problem, contains 4 binary and 14 real variables, 20 inequality and 2 equality constraints. 
 
Let $(z^*_{1|t},\ldots,z^*_{\tilde{p}|t})$ be the optimal solution computed for the binary variables
at time $t$, where $\tilde{p}=T\cdot \tilde{s}$. We exploit the shifted binary values optimized at the previous step
as the binary warm start, in particular $z_{1|t+1}=z^*_{\tilde{s}+1|t}$,
\ldots, $z_{\tilde{p}-\tilde{s}|t+1}=z^*_{\tilde{p}|t}$, by adopting Algorithm~\ref{algo:warmstart-sos1}.

We use the tolerances $ \epsilon_{V},~\epsilon_{G}=10^{-4}$, $\epsilon_{I}=10^{-2}$, and 
add the regularization term $10^{-2}I$, which makes the Hessian $Q$ positive definite. 
The average number of QPs solved and time is shown in Table~\ref{table:SYSID_results}, where \textsf{miqpGPAD} refers to  Algorithm~\ref{algo:MIQP-GPAD}, \textsf{miqpGPAD*} refers to Algorithms~\ref{algo:MIQP-GPAD} coupled with Algorithm~\ref{algo:warmstart-binary}  
and Algorithm~\ref{algo:warmstart-sos1} for warm-starting binary variables corresponding to the PWA model with $\tilde{s}=2$ modes. 
The average number of QP subproblems solved are reduced by almost 3 times, demonstrates the effectiveness of binary warm start coupled with~\textsf{noQP}, \textsf{noQP}$_{\TT_0}$, \textsf{noQP}$_{\TT_1}$ for modes of the PWA model under consideration.
\begin{table}[h!]
    \centerline{\begin{tabular}{ccc}
            \hline
            Solver             & Time (ms)& No. of subproblems \\ \hline
            GUROBI             & 2.0    & --\\
            \textsf{miqpGPAD}  & 5.7    & 13\\
            \textsf{miqpGPAD*} & 3.4    & 4\\ \hline
    \end{tabular}}
    \caption{Average CPU time and no. of QP subproblems solved for processing $N=1000$ training samples. }
    \label{table:SYSID_results}
\end{table}

\subsection{Heuristic approach --- Hybrid vehicle example}

The hybrid vehicle example from~\cite{TMBB17} consists of the combination of a battery, an electric motor/generator, and an engine. For a given power demand profile $P_{t}^{des}$, $t=0, \ldots,T-1$, the objective is to plan the battery power $P_{t}^{batt}$ and engine power $P_{t}^{eng}$ for the time interval $t=0, \ldots,T-1$ such that 
$P_{t}^{batt}+P_{t}^{eng} \geq P_{t}^{des}$.
Let $E_{t}$ be the energy of the battery at time $t$, $E_{t+1}=E_{t}-\tau P_{t}^{batt}$, where $\tau$ is the sampling time
and $0 \leq E_{t} \leq E^{max}$.  The fuel cost is given by $f(P_{t}^{eng},z_{t})$ where $f(P,z) =\alpha P^{2} + \beta P + \gamma z$, the constraint on power is $0 \leq P_{t}^{eng} \leq P^{max}z_{t}$. The optimal control problem to solve is the following:
\begin{eqnarray*}
	\begin{aligned}
		{\text{min}}
		& \quad \eta (E_{T}-E^{max})^{2}+\sum_{t=0}^{T-1}f(P_{t}^{eng},z_{t})\\&\quad +\delta \max(z_{t}-z_{t-1},0) \\
		\text{s.t.} & \quad E_{t+1}=E_{t}-\tau P_{t}^{batt}\\
		&   \quad P_{t}^{batt}+P_{t}^{eng} \geq P_{t}^{des}\\
		&  	\quad z_{t} \in \{0,1\},\ t=0, \ldots,T-1
	\end{aligned}
\end{eqnarray*}
where $ P_{t}^{batt} $, $P_{t}^{eng}$, $z_{t}$ (engine on/off) and $E_{t}$ are the optimization variables. The weight $\delta\geq0$ penalizes the engine going from the off to the on state.
By choosing $T = 72$ steps, the resulting MIQP problem has $n=862$ optimization variables, $p=72$ binary variables, $m=503$ inequality constraints, $q=575$ equality constraints.

The cost calculated by the ADMM-based heuristic approach of~\cite{TMBB17}, denoted as \textsf{miqpADMM}, is 138.1 and is obtained in 0.40~s for preconditioning, 3.55~s for solving the problem. The solver GUROBI computes the optimal cost $V^*=135.9$ 
in 21.05~s\footnote{\url{https://github.com/cvxgrp/miqp_admm/tree/master/matlab/vehicle.m}}. 

Algorithms~\ref{algo:GPAD}, \ref{algo:heuristic}, and~\ref{algo:Infeasibility detection}, heuristic approach
collectively denoted as \textsf{miqpGPAD-H}, are
implemented in interpreted MATLAB code for fairness of comparison with \textsf{miqpADMM}, 
and the resulting performance is reported in 
Table~\ref{tab:hybrid vehicle miqpgpad} for different values of 
the feasibility tolerance $\epsilon_{G}$ and optimality tolerance 
$\epsilon_{V}$.
Whenever condition~\eqref{eq:restart condition} is satisfied at a given iteration $k$, rather than restarting the values of $\beta_k$ we just assign $w_k\leftarrow \lambda_k$, $w_{eq,k}\leftarrow\nu_k$ in Step~\ref{eq:gpad steps w_k} (of Algorithm~\ref{algo:GPAD}) for that iteration. 
We add $10^{-3}I$ to the primal Hessian matrix $Q$ and
use $\epsilon_{I}=10^{-2}$ in Algorithm~\ref{algo:Infeasibility detection}.

\begin{figure}[ht!]
	\centering
	\includegraphics{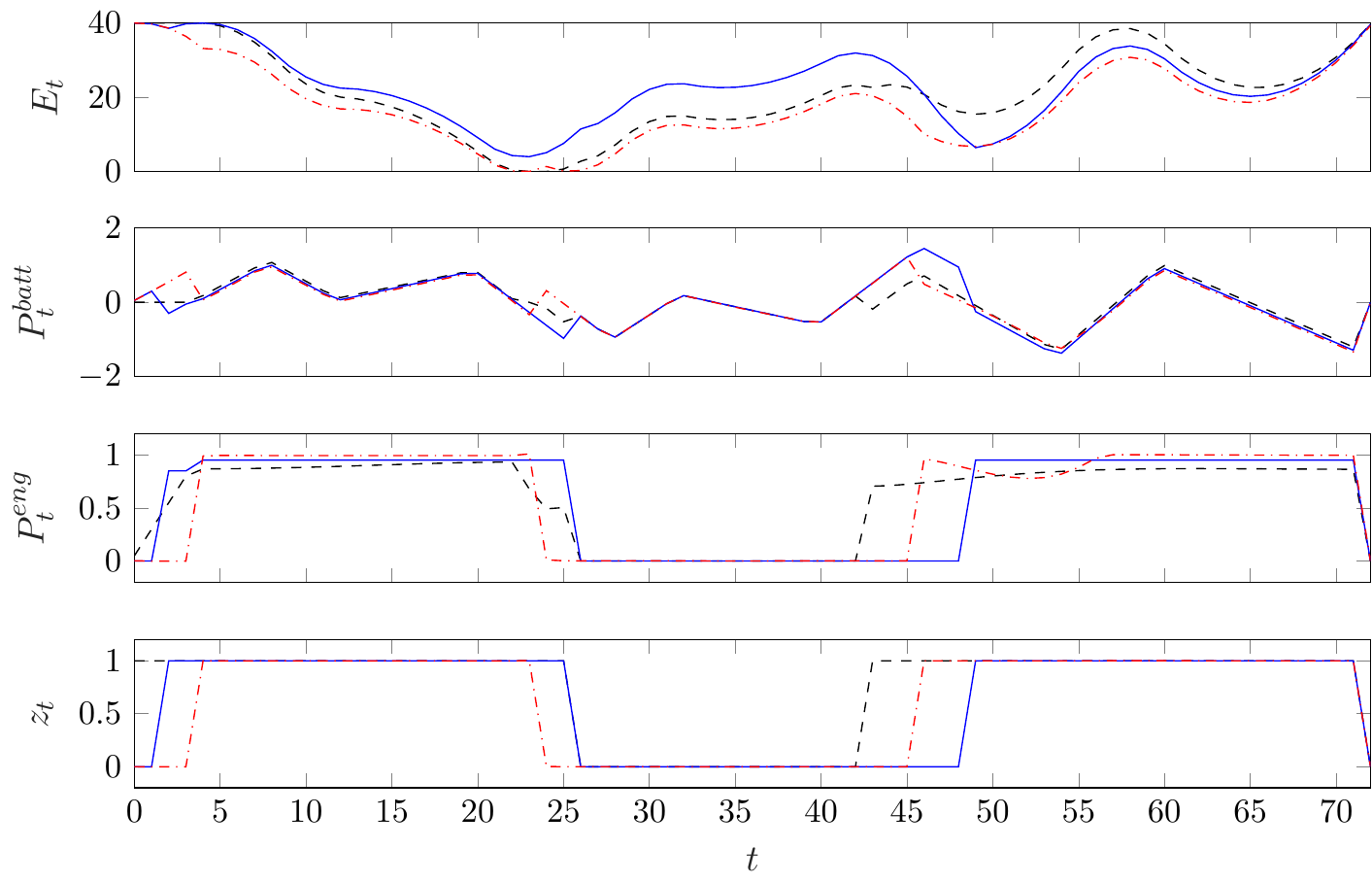}
	\caption{Battery energy, battery power, engine power and engine on/off signals versus time: GUROBI (solid blue line),~\textsf{miqpGPAD-H} (dash-dotted red line), and~\textsf{miqpADMM} (dashed black line).}
	\label{fig:vehicle_1_2}
\end{figure}

\begin{table}[h!]
	\centerline{\small\begin{tabular}{*4c}
			\hline\
			$ \epsilon_{V} $, $ \epsilon_{G} $ & Cost & Precond., Solving  &Constr. violation\\ \hline
			$10^{-2}$, $10^{-2}$& {131.5}  & {1.19, 3.81 s}  & $3.62 \cdot 10^{-1}$\\ \hline
			$10^{-2}$, $10^{-3}$& {135.9}  & {1.12, 8.62 s}  & $8.09 \cdot 10^{-3}$\\ \hline
			$10^{-3}$, $10^{-3}$& {135.9}  & {1.07, 8.67 s}  & $7.98 \cdot 10^{-3}$\\ \hline
			$10^{-3}$, $10^{-4}$& {136.0}  & {1.14, 13.37 s} & $2.59 \cdot 10^{-3}$\\ \hline
			$10^{-4}$, $10^{-3}$& {136.1}  & {1.10, 15.87 s} & $1.08 \cdot 10^{-4}$\\ \hline
		\end{tabular}
	}\vspace*{.5cm}
	\caption{Performance comparison with different values of $ \epsilon_{V} $, $ \epsilon_{G} $ for \textsf{miqpGPAD-H}.}
	\label{tab:hybrid vehicle miqpgpad}
\end{table}
Fig.~\ref{fig:vehicle_1_2} shows the trajectories obtained with \textsf{miqpGPAD-H} for $\epsilon_{V}=10^{-2}$, $\epsilon_{G}=10^{-3}$ and compare
them with the ones obtained by GUROBI and \textsf{miqpADMM}. 
It is apparent that the proposed heuristic approach
keeps the quality of the solution over a sufficient level for the practical application at hand, is computationally faster than GUROBI and comparable to \textsf{miqpADMM}, and very simple to implement in an embedded control platform.

\subsection{Mid-way approach on random MIQPs}

\begin{figure}[!ht]
	\centering
	\begin{minipage}{0.45\textwidth}
	 	\centering
		\includegraphics[width=1\hsize, height=6.7cm]{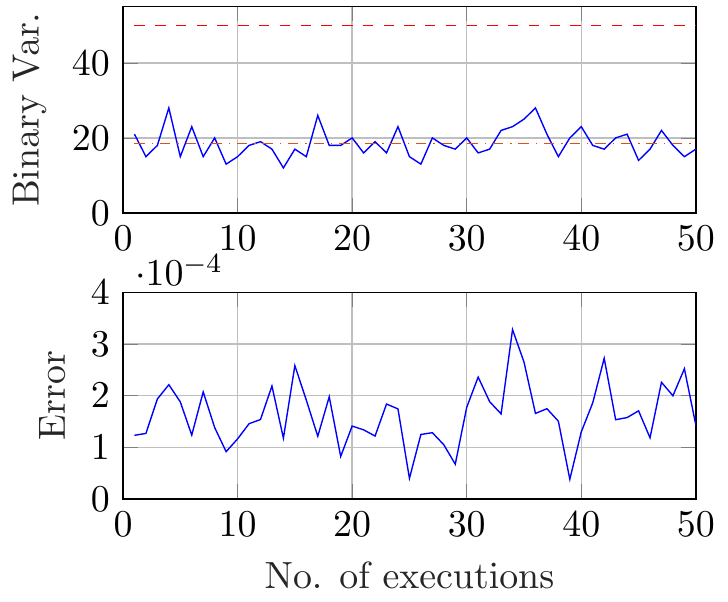}
	\end{minipage}
	\begin{minipage}{0.45\textwidth}
		\centering
		\includegraphics[width=1\hsize, height=6.7cm]{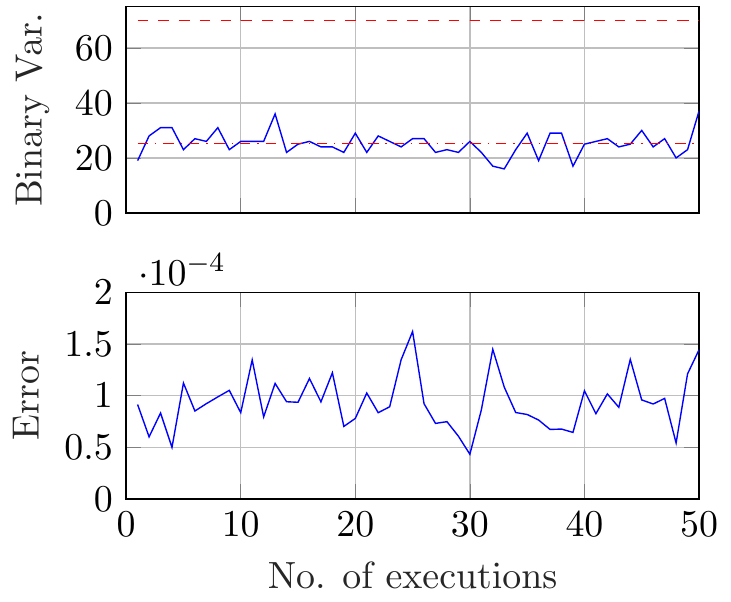}
	\end{minipage}
	\caption{ Performance of mid-way approach for considered two cases (first and second, left and right respectively), Top panel: Number of binary variables:~\textsf{miqpGPAD-H} (solid blue line), ~mean over $50$ runs (dash-dotted red line),~actual number of binary variables (dashed red line), Bottom panel: norm of error in decision variables.
	}
	\label{fig:heu_70bin_50runs_plot}
\end{figure}

We test the mid-way heuristic approach, denoted as \textsf{miqpGPAD-mH} on randomly generated MIQP problems
on two numerical experiments: first with $n=60$ variables, $m=60$ inequality constraints, $q=8$ equality constraints, and $p=50$ binary variables; and second with $n=90$ variables, $p=70$ binary variables, $m=80$ inequality and $q=8$ equality constraints with $ \epsilon_{V},\epsilon_{G}=1e-5,\epsilon_{I}=1e-2$.
The tolerances used in Algorithm~\ref{algo:heuristic} are $\epsilon_{\bar \ell} = 0.01$ and $\epsilon_{\bar u} = 0.99$, and the results for both experiments are shown in Fig.~\ref{fig:heu_70bin_50runs_plot}, through left and right panels respectively. In both cases, it is observed that over 50 runs, the mean value of actual number of binary variables solved using B\&B are 18 and 25 respectively, which is approximately 37\% of the original number of binaries.		

\subsection{Heuristic approaches on ARX model segmentation}

We apply the proposed heuristic approaches to an ARX model segmentation problem which can be solved using sparse optimization with $\ell_0$-norm penalties~\cite{Piga20133646}.

The problem of finding the sparsest solution satisfying a system of linear equations $A\theta=b$ with more unknowns than equalities can be given as
\begin{equation}
\min_{\theta} \|\theta\|_{0}~~\text{s.t.}~~A\theta=b
\label{eq: zero norm st equalities}
\end{equation}
where $A \in \mathbb{R}^{m \times n}$ with $m<n$, $b \in \mathbb{R}^m$, $\theta \in \mathbb{R}^n$ and $\|\cdot\|_{0}$ is the $\ell_0$-seminorm, that is the number of nonzero components of its argument. 
To take into account that the entries of matrix $A$ and/or the components of vector $b$ can be affected by noise,
Problem~\eqref{eq: zero norm st equalities} is relaxed to 
the Lagrangian form
\begin{equation}
\min_{\theta} \quad \|A\theta-b\|^{2}_{2}+\gamma\|\theta\|_{0} \label{eq: zero norm lagrangian form}
\end{equation}
where $\gamma>0$ is a tuning parameter. 

Problems~\eqref{eq: zero norm st equalities} 
and~\eqref{eq: zero norm lagrangian form} are nonconvex combinatorial problems. A common approximation is to replace the $\ell_{0}$-seminorm with the $\ell_{1}$-norm, which leads to solving a convex optimization problem. However, the exact solution of the original problem may not be achieved by such a relaxation. 
An exact solution of Problem~\eqref{eq: zero norm lagrangian form} can be retrieved by solving the following MIQP problem
\begin{subequations}
\begin{eqnarray}
&\min_{\theta, \omega} \quad \|A\theta-b\|^{2}_{2}+\gamma\sum_{i=1}^{n}\omega_{i} \label{eq:l0 as MIQP}\\
&\text{s.t.} \quad m \omega_{i} \leq \theta_{i} \leq M\omega_{i}~~i=1,\ldots,n
\end{eqnarray}\label{eq:l0_as_MIQP}
\end{subequations}
where $\omega_{i} \in \{0,1\}$ for $i=1,\ldots,n$, $M$ and $m$ are known upper and lower bounds on the components $\theta_{i}$.
It can be easily shown that Problems~\eqref{eq: zero norm lagrangian form} and~\eqref{eq:l0_as_MIQP} are equivalent\footnote{In the MIQP formulation, $\|\theta\|_{0}$ is replaced by $\sum_{i=1}^{n}\omega_{i}$,  and  the optimum $\left\lbrace\theta^{*}, \omega^{*}\right\rbrace$ is equal to that of Problem~\eqref{eq: zero norm lagrangian form}}.

The idea of using MIQP for sparse optimization using $\ell_{0}$-penalties can be applied to ARX model segmentation problems. Consider the ARX model
\begin{align}
y(t) = \varphi'(t)\theta+e(t)
\end{align}
where $\varphi(t)  = [-y(t-1)\ \ldots\ -y(t-n_{a})\ u(t-n_{k}-1)\ \ldots\ u(t-n_k-n_b)]'$ is the regressor vector and 
$\theta = [a_{1}\ \ldots\ a_{n_a}\ b_{1}\ \ldots\ b_{n_{b}}]'\in \rr^{n}$ is the vector of unknown parameters. A time-varying model can be equivalently represented as 
\begin{equation*}
y(t)= \varphi'(t)\theta(t)
\end{equation*}
In ARX model segmentation problems we assume that $\theta$ is time-varying
and piecewise constant,
\begin{equation*}
\theta(t) = \theta_{k}, \quad t_{k} < t \leq t_{k+1}
\end{equation*} 
We want to estimate $\theta(t)$ from a given noise-corrupted dataset of $N$ samples
$\left\lbrace y(t), u(t) \right\rbrace_{t=1}^{N}$.

We consider the example from \texttt{iddemo6m.mat} in the System Identification
Toolbox~\cite{ljung1988system, ohlsson_ljung20101107}, consisting of the system
\begin{align}\label{sys:transport_delay}
y(t)+0.9y(t-1)=u(t-n_{k})+e(t)
\end{align}
in which the value of  transport delay $n_{k}$ changes from 2 to 1 at $t=20$. 
The input $u$ is a pseudo-random binary sequence assuming values $\pm1$, $e(t)$ has variance $0.1$. 
The  ARX model  $y(t)= -ay(t-1)+b_{1}u(t-1)+b_{2}u(t-2)=\varphi(t)'\theta(t)$
is used to estimate  $\theta(t)=[a\ b_{1}\ b_{2}]'$, which equivalently 
represents~\eqref{sys:transport_delay} when $(b_1,b_2)$ switches between
$(1,0)$ and $(0,1)$.

The estimation problem is formulated as the following sparse-optimization problem with $\ell_{0}$-penalty
\begin{equation}
\min_{\theta(t)} \quad \sum_{t=1}^{N}(y(t)-\varphi(t)'\theta(t))^2 +\gamma\sum_{t=2}^{N}\|\theta(t)-\theta(t-1)\|_{0} \label{l0 norm problem-delta theta}
\end{equation}  

The term $\|\theta(t)-\theta(t-1)\|_{0}$ penalizes changes 
of model parameters over time. The hyper-parameter $\gamma$ can be tuned to achieve a trade-off between model fit and frequency of variation of model parameters.
Problem~\eqref{l0 norm problem-delta theta} can expressed as an MIQP problem~\eqref{eq:l0_as_MIQP}.

The obtained MIQP problem consists of $n=240$ variables, out of which $p=120$ are binary variables, and $q=240$ inequality constraints. 
The regularization term $10^{-2}I$ is added to the Hessian matrix and $M=-m=1$.
This problem is solved using the heuristic \textsf{miqpGPAD-H} and the mid-way heuristic \textsf{miqpGPAD-mH} approaches. 

\begin{table}[ht!]
	\centerline{\begin{tabular}{*4l}
			\hline
			Solver                  &$\gamma$ & Cost  & Time (s) \\ \hline
			GUROBI                  &0.5 & 46.17 & 0.70     \\
			\textsf{miqpGPAD-H}     & & 47.15 & 0.47     \\ \hline
			GUROBI                  &0.7 & 46.97 & 0.17     \\
			\textsf{miqpGPAD-H}     & & 47.95 & 0.89     \\
			\hline
	\end{tabular}}
	\caption{Performance comparison for ARX model segmentation problem using heuristic approach without using B\&B.}
	\label{table:ARX_seg_per_comp}
\end{table}
\begin{figure}[ht!]
	\centering
	\includegraphics[width=0.78\hsize, height=8.4cm]{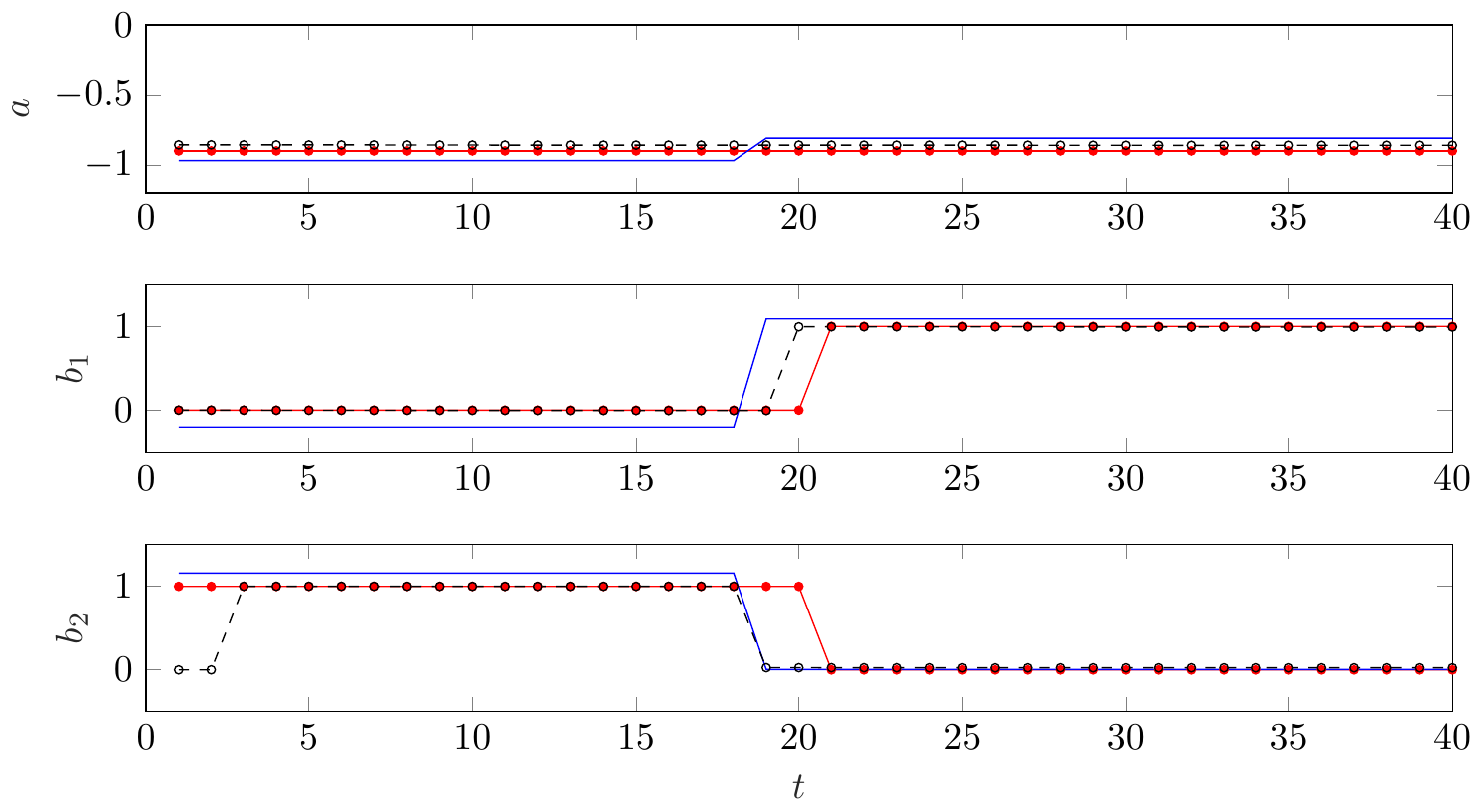}
	\caption{Performance comparison for ARX model segmentation problem with $\gamma=0.7$, $\epsilon=0.2$ :~true value (dotted solid red line),~\texttt{segment} (solid blue line),~\textsf{miqpGPAD-H} (dashed-dotted black line).}	
	\label{fig:plot_ARX_seg_heu}
\end{figure}
The results using \textsf{miqpGPAD-H} with $\epsilon_{V}=10^{-5},\epsilon_{G}=10^{-2},~\epsilon_{I}=10^{-2}$
are summarized in Table~\ref{table:ARX_seg_per_comp}.
For \textsf{miqpGPAD-mH}, the experiments are carried out with varying values of $\epsilon_{\bar{\ell}}=1-\epsilon_{\bar{u}}=\epsilon$ and the results are shown in Table~\ref{table:ARX_seg_per_comp_mheu}, which also reports the number $p$ of binary variables 
used in B\&B. The tolerance for the B\&B algorithm \textsf{miqpGPAD} are $ \epsilon_{V},\epsilon_{G}=10^{-5},~\epsilon_{I}=10^{-2}$. 
\begin{table}[ht!]
	\centerline{\begin{tabular}{*6l}
		\hline
		Solver               & $\gamma$ & $\epsilon$ & $p$ &  Cost  & Time (s) \\ \hline
		\textsf{miqpGPAD-mH} & 0.5      & 0.3        & 5   & 46.39 & 0.47, 0.70 \\
							 &          & 0.2        & 8   & 46.17 & 0.47, 1.27  \\
							 &          & 0.1        & 12  & 46.17 & 0.47, 1.85  \\
		\hline
		\textsf{miqpGPAD-mH} & 0.7       & 0.3        & 5   &  47.19 & 0.89, 0.58 \\
							 &          & 0.2        & 7   &  46.97 & 0.89, 0.86  \\
				   	     	 &          & 0.1        & 10  &  46.97 & 0.89, 1.42  \\ 
		\hline
	\end{tabular}}
	\caption{Performance comparison for ARX model segmentation problem using the mid-way heuristic approach, $p$ denotes the number of binary variables against total of 120.}
	\label{table:ARX_seg_per_comp_mheu}
\end{table}
\begin{figure}[ht!]
	\centering
	\includegraphics[width=0.78\hsize, height=8.4cm]{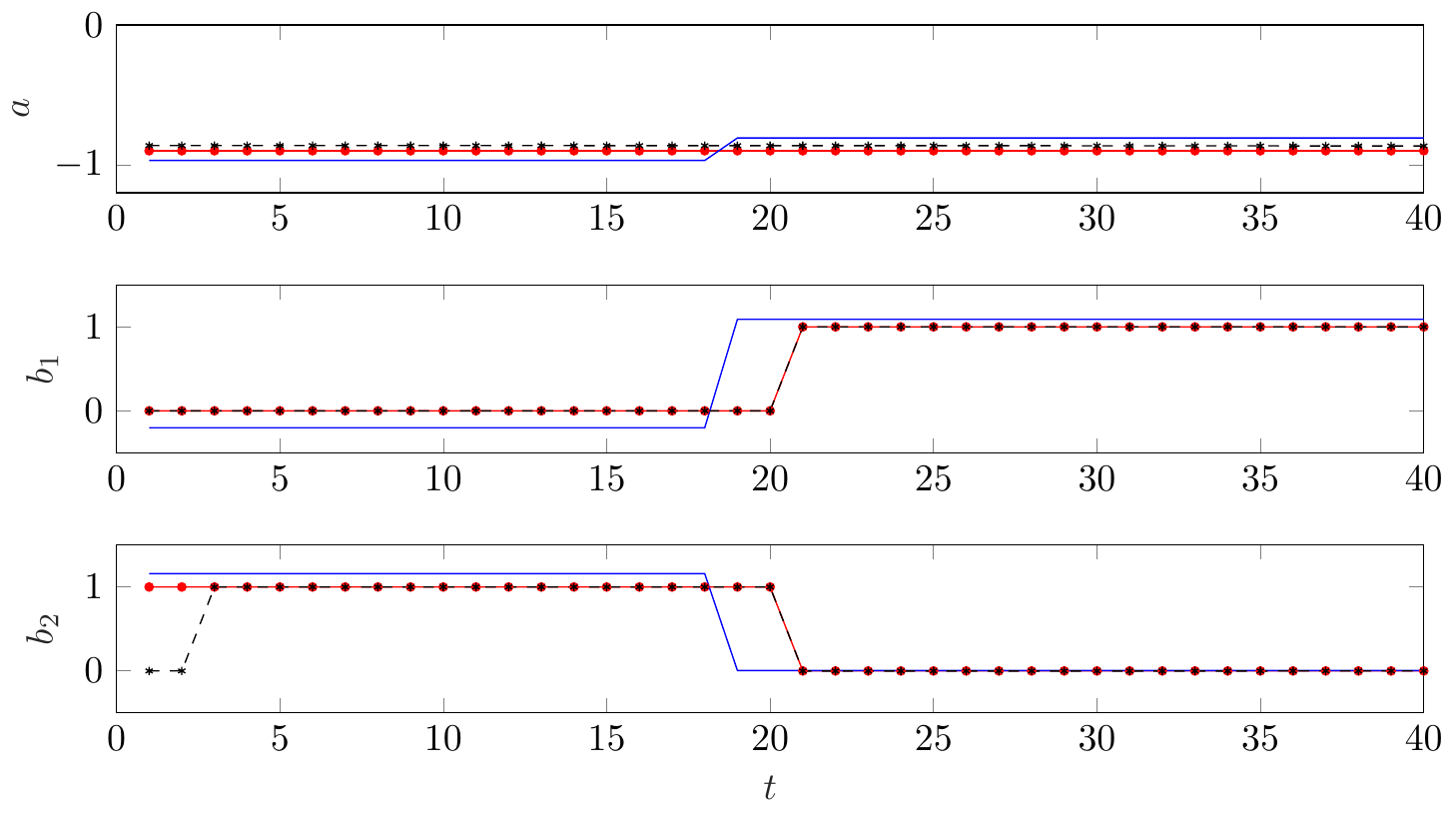}
	\caption{Performance comparison for ARX model segmentation problem with $\gamma=0.7$, $\epsilon=0.2$ :~true value (dotted solid red line),~\texttt{segment} (solid blue line),~\textsf{miqpGPAD-mH} (dashed-asterisk black line).}	
	\label{fig:plot_ARX_seg_mheu}
\end{figure}

The estimated parameter values $a$, $b_{1}$ and $b_{2}$ are compared with the true values and
the values produced by \texttt{segment} function~\cite{ljung1988system}, as shown in Figures~\ref{fig:plot_ARX_seg_heu} and~\ref{fig:plot_ARX_seg_mheu} for \textsf{miqpGPAD-H} and \textsf{miqpGPAD-mH} respectively.

The~\textsf{miqpGPAD-H} approach produces a good suboptimal solution within a short time using a very simple piece of code. 
Moreover, a further improvement over the solution of~\textsf{miqpGPAD-H} is demonstrated using~\textsf{miqpGPAD-mH}, which employs B\&B for solving a few binary variables that comes with trading-off more computation efforts and time. As~\textsf{miqpGPAD-mH} explores how to solve significantly fewer binaries using B\&B, this improvement is achieved with just approximately 7\% of the total number of binary variables to branch on, and produces a solution which is further close to the optimal solution within the specified tolerance values.

\section{Conclusions}

In this paper we have presented very simple and practical, exact and heuristic approaches to solve MIQP problems based on accelerated gradient projection applied on the dual QP relaxations. Furthermore, we have proposed warm-starting strategies, both in solving QP relaxations and in prioritizing the combination of binary constraints during B\&B. A comprehensive study of warm-starting binary variables framework combined with B\&B have been presented,  which formally supports any kind of presolving techniques as well as MIQP with binary variables correspond to an exclusive-or condition (SOS1 type constraints). It inherently prioritize the tree as per the warm-start provided, avoids solving their parents while exploring the remaining nodes, without compromising on the optimality of the solution.
In spite of their simplicity,
the proposed approaches were shown quite effective
in addressing relatively small-size MIQP problems that arise from embedded control and estimation applications. The heuristic method without B\&B can often provide solutions close to optimality with much reduced coding and computation efforts, which can be also used to provide a good initial upper-bound on the optimal
cost, in case B\&B is executed to determine the optimal solution. By combining binary warm-starting with the heuristic approaches, the number of solved QP relaxations can be 
reduced. 
Current research is devoted to combining warm-starting ideas presented in this paper with prediction techniques
based on artificial neural networks~\cite{MB2019_learning_ws}.

\section*{Acknowledgement}
	The authors would like to thank Dr. Dario Piga and Dr. Manas Mejari for their help with the ARX model segmentation example.

\bibliographystyle{plain}       
\bibliography{miqpgpad_extn_arxiv}  

\end{document}